\documentclass[twoside,11pt]{article}
\usepackage{amsmath,amsthm}
\let\amsproof\proof
\let\endamsproof\endproof
\usepackage[preprint,nohyperref,abbrvbib]{jmlr2e}
\let\proof\amsproof
\let\endproof\endamsproof
\usepackage[T1]{fontenc}
\usepackage{lmodern}
\usepackage{mathtools,booktabs,microtype,mathrsfs}
\usepackage[caption=false]{subfig}
\usepackage{placeins}
\usepackage{xcolor}
\usepackage[hidelinks]{hyperref}
\makeatletter\g@addto@macro{\UrlBreaks}{\do\-}\makeatother
\hypersetup{pdftitle={Neural Operator Surrogates for Landweber Iteration in Solving Inverse Problems},pdfauthor={Kuo Gai, Marco A. Iglesias, Shuai Lu, Frank Werner},pdfsubject={Neural operators, nonlinear inverse problems, iterative regularisation, Landweber iteration, derivative-informed learning}}
\ShortHeadings{Neural Operator Surrogates for Landweber Iteration}{Gai, Iglesias, Lu and Werner}
\firstpageno{1}
\newcommand{\M}{\mathcal M}
\newcommand{\Y}{\mathcal Y}
\newcommand{\F}{\mathcal F}
\newcommand{\Adm}{\mathcal K}

\newcommand{\Lop}{\mathcal L}
\newcommand{\epsop}{\varepsilon_{\mathrm{op}}}

\newcommand{\Deff}{\Delta}
\newcommand{\ip}[2]{\langle #1,#2\rangle}

\DeclareMathOperator{\ran}{ran}
\DeclareMathOperator{\dom}{dom}
\DeclareMathOperator{\dist}{dist}
\DeclareMathOperator{\Lip}{Lip}
\DeclareMathOperator*{\argmin}{argmin}
\numberwithin{theorem}{section}
\theoremstyle{definition}
\newtheorem{assumption}[theorem]{Assumption}
\title{Neural Operator Surrogates for Landweber Iteration in Solving Inverse Problems}
\author{\name Kuo Gai \email 25110180015@m.fudan.edu.cn\\
       \addr School of Mathematical Sciences, SKLCAM and LMNS\\
       Fudan University, Shanghai, China
       \AND
       \name Marco A. Iglesias \email marco.iglesias@nottingham.ac.uk\\
       \addr School of Mathematical Sciences\\
       University of Nottingham, Nottingham, UK
       \AND
       \name Shuai Lu \email slu@fudan.edu.cn\\
       \addr School of Mathematical Sciences, SKLCAM and LMNS\\
       Fudan University, Shanghai, China
       \AND
       \name Frank Werner \email frank.werner@uni-wuerzburg.de\\
       \addr Institute of Mathematics\\
       University of W\"urzburg, W\"urzburg, Germany}
\begin{document}
\maketitle
\thispagestyle{plain}
\begin{abstract}
Neural operators can accurately approximate well-posed mappings, but for unbounded or discontinuous operators, common in inverse problems, we verify that their network parameters often blow up in architectures such as DeepONet and FNO. To solve inverse problems stably, we use Landweber iteration, which requires only a forward map and the adjoint of its derivative (linearisation). To reduce repeated PDE forward evaluations, we replace the forward operator and its derivative with a differentiable neural-operator surrogate. We show that encoder–network–decoder operators, including generalised DeepONets and fully coupled Fourier operators, approximate both the forward map and its derivative. A modified discrepancy principle ensures decreasing error before finite stopping. With a compact preconditioner, the iterates remain in a compact set where accurate surrogates exist, and the stopped reconstructions converge as data and surrogate errors vanish under a regularity condition. We also relate the theory to derivative-informed training along finite directions, showing that the preconditioner determines which derivative errors affect the Landweber iteration update. Experiments on identifying a reaction coefficient in an elliptic equation show that smoothing reduces late-iteration deterioration for less accurate surrogates but slows sharp-feature recovery, and that smaller derivative errors along training directions do not by themselves improve reconstruction accuracy.

\end{abstract}
\section{Introduction}
\label{sec:intro}
\label{sec:framework}

Many computational tasks involving parametric partial differential equations
(PDEs) require repeated evaluation of the same parameter-to-solution or
parameter-to-observable map. Examples include uncertainty quantification,
PDE-constrained optimisation and control, experimental design and parameter
identification. When each evaluation requires the solution of a large-scale
PDE, the cost can become prohibitive, particularly when sensitivities of the
forward map are needed as well as its values.

Neural operators offer one way of amortising this cost. A neural operator is
trained offline to approximate a map between function spaces and is then
evaluated at comparatively low online cost
\citep{Lu2021,Li2021,Kovachki2023}. Learned operators have been embedded in
PDE-constrained optimisation and control \citep{Hwang2022}, Bayesian
experimental design \citep{GoChenDINO} and Bayesian parameter identification,
where posterior exploration requires many forward and derivative evaluations
\citep{DINOMCMC2025}. Derivative-informed surrogates developed for such
many-query settings learn the forward map together with its derivative
\citep{OLearyRoseberry2022,DINO2024}. Neural-operator surrogates have also been
used in ensemble-based inversion of moving-boundary problems
\citep{IglesiasDeepONet2026}.

Neural operators can also be trained to map observations directly to the
unknown parameters. Neural Inverse Operators, introduced by
\citet{Molinaro2023}, learn operator-to-function maps for PDE inverse problems
with operator-valued measurements. For the Calder\'on problem,
\citet{deHoop2025} construct a noise-aware neural-operator approximation using
an extension of the inverse conductivity map that preserves its stability
properties. In inverse imaging, \citet{Auras2026} interpret U-Nets, which are
widely used to post-process a crude reconstruction such as the backprojected
data, as neural operators, and study how such networks behave when the
discretisation is refined, with limited-angle computed tomography as the
application. These approaches rely on structure, or on a regularised input,
that keeps the learned inverse map stable on the class of interest. Without such structure, the inverse map of an
ill-posed problem is unstable, and this has a direct consequence for neural
approximation (Proposition~\ref{prop:inverse-blowup} in
Section~\ref{app:inverse-blowup}). No continuous network can approximate a
discontinuous inverse map uniformly near a point of discontinuity to within
half the size of the jump, and networks that are accurate on training data
resolving the discontinuity must have Lipschitz constants that grow without
bound. For Fourier neural operators of fixed depth, some weights must then
diverge. In this work we therefore learn the forward map, which in PDE
parameter identification is typically continuously differentiable, and
control the instability of the inversion by regularisation.

On the other hand, the use of learned forward maps within \emph{classical regularisation of
ill-posed inverse problems} has received less attention. In this setting the
learned operator is part of the regularisation method itself, and replacing
the forward map by an approximation may alter stability, convergence and the
criterion by which the iteration should be stopped. \citet{ScherzerVuYan2025}
study Tikhonov regularisation with an approximate neural forward operator and
show how observational and operator errors must be balanced with the
regularisation parameter. Landweber-type iterations with perturbed or noisy
forward operators have also been analysed \citep[e.g.][]{GuDongWang2025}, and
learned components have been combined with Landweber iteration
\citep{Aspri2020}. To our knowledge, however, the use of neural-operator
surrogates within iterative regularisation, where the derivatives that
determine the search directions are also approximate, remains largely
unexplored.

Iterative regularisation methods evaluate the forward map and, for
derivative-based schemes, its linearisation or adjoint at every step. When
these operations are expensive, replacing them by a differentiable learned
operator could reduce the online cost substantially once the offline training
cost has been amortised. This is most relevant when many related inverse
problems must be solved, or when reconstructions are needed quickly, for
example in real-time monitoring. We study this question for nonlinear
Landweber iteration. Landweber is among the simplest derivative-based
iterative regularisation methods and is typically slower than Newton,
Gauss--Newton or Levenberg--Marquardt-type schemes. Inexpensive surrogate
evaluations change the balance between iteration count and cost per
iteration, which makes this slower method more attractive. Landweber iteration also
provides a transparent setting in which the effect of approximating both the
forward map and its derivative can be analysed. The analysis is formulated for
nonlinear maps between Hilbert spaces, and PDE parameter identification
provides the computational motivation and the numerical examples.

Let $\M$ and $\Y$ be real separable Hilbert spaces and let $\F:\operatorname{dom}(\F)\subseteq\M\longrightarrow\Y$
denote the forward operator. Given exact data $y^\dagger=\F(m^\dagger)$ and
noisy observations $y^\delta$ satisfying
\begin{equation}
\|y^\delta-y^\dagger\|_{\Y}\leq\delta,
\label{eq:observations}
\end{equation}
the nonlinear Landweber iteration takes the form
\begin{equation}
m_{k+1}
=
m_k-\omega D\F(m_k)^*
\bigl(\F(m_k)-y^\delta\bigr),
\label{eq:exact-landweber}
\end{equation}
together with a stopping rule. Under standard local assumptions, including a
tangential cone condition, Landweber iteration stopped by the discrepancy
principle is a regularisation method for nonlinear ill-posed problems
\citep{HankeNeubauerScherzer1995,EHN1996,KNS2008}. As the iteration formula \eqref{eq:exact-landweber} shows, the Landweber iteration is an \emph{explicit regularisation method} that involves only the forward mapping and its linearisation, both of which can be well approximated by suitable learned neural operators. Replacing $\F$ by a
differentiable learned operator $\F_\theta$, and differentiating the same
surrogate, gives
\begin{equation}
m_{k+1}
=
m_k-\omega D\F_\theta(m_k)^*
\bigl(\F_\theta(m_k)-y^\delta\bigr).
\label{eq:basic-landweber}
\end{equation}
The central question of this paper is whether this substitution preserves the
regularising properties of the original method.

Two approximation errors enter the iteration. Error in $\F_\theta$ perturbs
the residual against which the data are fitted, and error in $D\F_\theta$
perturbs the direction in which the parameter is updated. Accurate forward
prediction alone is therefore not sufficient, because a surrogate with
inaccurate derivatives may produce poor Landweber directions. The usual
discrepancy principle also cannot be applied unchanged, since a small residual
for $\F_\theta$ need not imply that the physical forward equation is satisfied
to the same accuracy.

Our contributions are as follows.
\begin{itemize}
 \item We show that near a discontinuity, no continuous neural network can uniformly approximate a discontinuous inverse map, which renders neural operator surrogates for the inverse map infeasible (Section~\ref{app:inverse-blowup}). Meanwhile, we identify encoder--network--decoder neural operators that
 approximate a nonlinear forward map and its derivative simultaneously,
 first along compact sets of directions and then, for compact derivatives,
 in operator norm (Section~\ref{sec:approx}). Generalised DeepONets and fully
 coupled Fourier operators realise this class exactly
 (Section~\ref{sec:architectures}).
 \item For a surrogate with separate forward and derivative error bounds, we
 introduce an effective discrepancy level that accounts for observational,
 forward and derivative errors, and prove descent before stopping and finite
 termination of the resulting discrepancy principle
 (Section~\ref{sec:landweber}).
 \item We prove convergence of the stopped reconstructions as the data and
 surrogate errors vanish for the preconditioned iteration
 \[
  m_{k+1}=m_k-\omega\,C\,D\F_\theta(m_k)^*\bigl(\F_\theta(m_k)-y^\delta\bigr)
 \]
 with a compact operator $C$ (Theorem~\ref{thm:convergence}). Closed
 balls of the Hilbert space $\M_C$ induced by $C$ are compact in $\M$, descent
 in the $\M_C$ norm keeps the iterates in such a ball, and the approximation
 theory supplies surrogates that are accurate on it. The price is the
 regularity condition $m^\dagger-m_0\in\M_C$. Without a compact
 preconditioner, the descent argument still applies wherever the error bounds
 hold, but the iterates are no longer confined to a compact set and
 convergence as the errors vanish is not established.
 \item We relate these requirements to practical derivative-informed training
 (Section~\ref{sec:training-preconditioning}). Supervision along a finite set
 of directions controls a finite-rank part of the derivative error, whereas
 the preconditioned iteration depends on the weighted error
 $(D\F_\theta-D\F)C^{1/2}$. Smoothing can attenuate poorly learned components
 of the search direction, but it also suppresses genuine fine-scale features
 and restricts the truths covered by the theory.
 \item We study these effects numerically for the identification of a
 reaction coefficient in an elliptic PDE (Section~\ref{sec:numerical-interpretation}),
 comparing standard and generalised architectures, training with and
 without the derivative loss, and different parameter-space geometries. The
 networks are trained with derivatives along a fixed set of empirical
 Karhunen--Lo\`eve (KL) directions, which is one practical choice among
 several covariance-, sensitivity- or random-projection-based options
 \citep{DINO2024,Qiu2024,Yao2025,Luo2026}. Transient Darcy-flow inversion with sparse space--time data
 (Section~\ref{sec:darcy}) illustrates the approach in a time-dependent
 setting with a standard architecture.
\end{itemize}
Conclusions and directions for further work are given in
Section~\ref{sec:conclusions}.
The appendices contain the approximation proofs, the standard Fourier Neural Operator (FNO) calculus and the proofs for the preconditioned geometry.

\section{Learning a discontinuous inverse map}
\label{app:inverse-blowup}

This section formalizes the statement made in the introduction that directly learning a discontinuous inverse map is unstable, thereby establishing the necessity of our proposed approach, which uses a neural-network-surrogated Landweber iteration to solve inverse problems stably.

Let $\mathcal R:\dom(\mathcal R)\subseteq\Y\to\M$ be an inversion or reconstruction map, for example, $\mathcal R=\F^{-1}$ on the range of an injective forward operator. The following is a formal definition of such a discontinuous operator. 
\begin{definition}[Discontinuous map]
\label{def:discontinuous-inverse}
The map $\mathcal R$ is discontinuous at $y_0\in\dom(\mathcal R)$ if there exist a constant $C>0$ and a sequence $\{y_n\}\subset\dom(\mathcal R)$ such that $y_n\to y_0$ in $\Y$, $y_n\ne y_0$ for all $n$, and
\[
\|\mathcal R(y_n)-\mathcal R(y_0)\|_\M\ge C \quad \text{for all } n.
\]
\end{definition}

Consequently, the following proposition shows that no continuous map can accurately approximate a discontinuous map, including Lipschitz continuous maps and generic neural operators. For simplicity, we prove the result only for FNO, but the argument holds for other architectures such as DeepONet. 
\begin{proposition}[Weight Blow-Up in Approximation of Discontinuous Maps]
\label{prop:inverse-blowup}
Let $\mathcal R$ be discontinuous at $y_0$, with $C$ and $\{y_n\}$ as in
Definition~\ref{def:discontinuous-inverse}. Then the following hold. 
\begin{enumerate}
 \item[(i)] Let $G:\dom(\mathcal R)\to\M$ be continuous at $y_0$. If
 $\|G(y)-\mathcal R(y)\|_\M\le\epsilon$ for $y=y_0$ and for all $y_n$,
 then there holds $\epsilon\ge C/2 > 0$.
 \item[(ii)] Let $0\le\epsilon_1 <C/2$, and for each $n$ let
 $G_n:\dom(\mathcal R)\to\M$ be Lipschitz continuous with
 $\|G_n(y)-\mathcal R(y)\|_\M\le\epsilon_1 $ for $y\in\{y_0,y_n\}$. Then
 \[
  \Lip(G_n)\ \ge\ \frac{C-2\epsilon_1 }{\|y_n-y_0\|_\Y}\ \longrightarrow\ \infty .
 \]
 \item[(iii)] 
 Suppose, in (ii), that each $G_n$ is a standard FNO with fixed depth and a globally Lipschitz activation, and that its lifting and projection networks have fixed depth. Let $\Y$ and $\M$ carry the $L^2$ or grid norms in which
 \eqref{eq:fno-lipschitz} in Appendix \ref{app:fno-calculus} holds. Then at least one of the weight-matrix or
 Fourier-multiplier norms in \eqref{eq:fno-lipschitz} is unbounded as
 $n\to\infty$.
\end{enumerate}
\end{proposition}
\begin{proof}
The first two arguments are well known in inverse problems, and we briefly outline the proof below. For any map $G$ and any $n$, the following inequality holds:
\[
 C\le\|\mathcal R(y_n)-\mathcal R(y_0)\|_\M
 \le\|\mathcal R(y_n)-G(y_n)\|_\M+\|G(y_n)-G(y_0)\|_\M
   +\|G(y_0)-\mathcal R(y_0)\|_\M .
\]

Then we can prove:
(i) the right-hand side is at most $2\epsilon+\|G(y_n)-G(y_0)\|_\M$, which
tends to $2\epsilon$ as $n\to\infty$ by continuity of $G$ at $y_0$. Hence $\epsilon\ge C/2>0$.

(ii) Choosing $G=G_n$, we see that the right-hand side of the above inequality is at most $2\epsilon_1+\Lip(G_n)\|y_n-y_0\|_\Y$. Rearranging yields a lower bound on $\Lip(G_n)$ whose right-hand side diverges as $y_n\to y_0$, since $C-2\epsilon_1>0$. 

(iii) By deriving the Lipschitz constant for FNO, as shown in the right-hand side of \eqref{eq:fno-lipschitz} in Appendix \ref{app:fno-calculus}, it is a product of the
listed norms and fixed constants. If all these norms were bounded in $n$, so
would be $\Lip(G_n)$, which contradicts the result of (ii).
\end{proof}

As we argued in the introduction and above, a neural operator cannot accurately approximate a discontinuous map, which in inverse problems is typically the inversion map. Nevertheless, we prove in the next section that if the forward map is well-posed, then generic neural operators can approximate the forward map and its linearisation to a certain accuracy.

\section{Simultaneous approximation of the forward map and its derivative}
\label{sec:approx}

Our proposed surrogate Landweber iteration \eqref{eq:basic-landweber} uses both the forward map and the adjoint
of its derivative, so we need neural operators that approximate both at once.
We work with a reduced-basis class built from a bounded linear encoder, a
finite-dimensional neural network and a bounded linear decoder. We first
approximate derivative actions on compact sets of directions, and then, under
additional assumptions, the derivative in operator norm, which is the bound
used in Section~\ref{sec:landweber}.

In this section, we allow $\M$ and $\Y$ to be Banach spaces.  The error-bound analysis of the surrogate Landweber iteration \eqref{eq:basic-landweber} in Section~\ref{sec:landweber} is carried out in Hilbert spaces, but it can be extended to Banach spaces by adapting the convergence analysis of \citet{SKHK2012}. An operator in the
reduced-basis class has the form
\begin{equation}
 \F_\theta(m)=\Phi f_\theta(\Psi m),\qquad
 \Psi\in\Lop(\M,\mathbb R^r),\quad
 \Phi\in\Lop(\mathbb R^s,\Y).
 \label{eq:rb}
\end{equation}
Here, the bounded linear encoder $\Psi$ maps the input $m$ to $r$
coordinates, the finite neural network
$f_\theta:\mathbb R^r\to\mathbb R^s$ maps these to $s$ output coefficients,
and the bounded linear decoder $\Phi$ reconstructs an element of $\Y$
from those coefficients. In the approximation results, the encoder, the decoder and their dimensions
are chosen together with the network, so greater accuracy may require richer
input and output representations as well as a larger network.

For the finite-dimensional network $f_\theta$, we use a nonpolynomial
activation $\sigma\in C^1(\mathbb R)$, such as GELU or $\tanh$.
Networks with this activation can approximate a finite-dimensional $C^1$
map and its first derivative uniformly on compact sets
\citep{Hornik1990,Pinkus1999}. This is the finite-dimensional step in the constructions below, and
Lemma~\ref{lem:finite-c1} in Appendix~\ref{app:approx-proofs} states the
version used here.

\subsection{Compact direction sets}

\citet[Theorem~4.4(ii)]{deFeo2026} shows that encoder--decoder neural
operators, a class that includes DeepONets, approximate a globally defined operator in his class $C_B^k$ and its first
$k$ derivatives simultaneously on compact sets
of inputs and directions. Theorem~\ref{thm:compact-directions} gives a local,
first-derivative version of this result for forward maps defined on an open
subset of the input space.
Appendix~\ref{app:approx-proofs} provides a self-contained proof and the
corresponding encoder--network--decoder construction.

A Banach space has the approximation property if its identity can be
approximated uniformly on every compact set by bounded finite-rank maps. In
the theorem below, this property provides finite-dimensional representations
of the inputs and directions in $\M$ and of the corresponding values and
derivative outputs in $\Y$.

For standard FNOs, simultaneous approximation of an
operator and its derivative is established by \citet[Theorem~3.2]{Yao2025}.
Section~\ref{sec:fno-comparison} compares that result with the ones below.

\begin{samepage}
\begin{theorem}[Approximation on compact inputs and directions]
\label{thm:compact-directions}
Suppose $\M$ and $\Y$ have the approximation property as defined in Definition \ref{def:ap}, $U\subset\M$ is open,
and $\F\in C^1(U,\Y)$. Let $\Adm\subset U$ and $H\subset\M$ be compact.
For every $\epsilon>0$, there exists an operator of the form
\eqref{eq:rb} such that
\begin{align}
 \sup_{m\in\Adm}\|\F(m)-\F_\theta(m)\|_\Y&<\epsilon,\label{eq:compact-value}\\
 \sup_{\substack{m\in\Adm\\h\in H}}
 \|D\F(m)h-D\F_\theta(m)h\|_\Y&<\epsilon.
 \label{eq:compact-direction}
\end{align}
\end{theorem}
\begin{proof}
See Appendix~\ref{proof:compact-directions}. 
\end{proof}
\end{samepage}

\subsection{Approximation of derivatives in operator norm}

Theorem~\ref{thm:compact-directions} controls the derivative error along a
prescribed compact set of directions. The Landweber analysis in
Section~\ref{sec:landweber} requires a bound over all directions of norm at
most one, that is, in the operator norm
\[
 \|D\F(m)-D\F_\theta(m)\|_{\Lop(\M,\Y)}
 =
 \sup_{\|h\|_\M\le1}
 \|(D\F(m)-D\F_\theta(m))h\|_\Y.
\]
Because the unit ball is not compact in infinite dimensions,
Theorem~\ref{thm:compact-directions} does not give this bound. The following
theorem does, under additional assumptions that include compactness of the
derivative.

\begin{theorem}[Uniform operator-norm approximation]
\label{thm:operator-approx}
Let $\M$ be a reflexive Banach space with a Schauder basis and let $\Y$ have
a Schauder basis. Suppose $U\subset\M$ is open,
$\F\in C^1(U,\Y)$, and $\Adm\subset U$ is compact. If $D\F(m)$ is compact
for every $m\in\Adm$, then for every $\epsilon>0$ there is an operator
of the form \eqref{eq:rb} such that
\[
 \sup_{m\in\Adm}\|\F(m)-\F_\theta(m)\|_\Y<\epsilon,\qquad
 \sup_{m\in\Adm}\|D\F(m)-D\F_\theta(m)\|_{\Lop(\M,\Y)}<\epsilon.
\]
\end{theorem}
\begin{proof}
See Appendix~\ref{proof:operator-approx}. 
\end{proof}

The compactness assumption in Theorem~\ref{thm:operator-approx} cannot be
removed for the class \eqref{eq:rb}. Its derivatives have the form $D\F_\theta(m)=\Phi\,Df_\theta(\Psi m)\,\Psi$
and therefore have finite rank, since their range lies in the range of
$\Phi$. An operator-norm limit of finite-rank operators is compact, so
approximating $D\F(m)$ arbitrarily well in operator norm requires $D\F(m)$ to
be compact. More generally, \citet[Theorem~4.2]{deFeo2026} shows that encoder--decoder architectures are not dense in $C^k$ with derivatives measured in operator norm whenever some bounded operator from the input to the output space cannot be approximated by finite-rank operators.

\begin{corollary}[Approximation of weighted derivatives]
\label{cor:weighted-approximation}
Suppose the assumptions of Theorem~\ref{thm:compact-directions} hold on a
compact input set $\Adm\subset U$, and let $T\in\Lop(\M)$ be compact. Then,
for every $\epsilon>0$, there exists an operator of the form \eqref{eq:rb}
such that
\[
 \sup_{m\in\Adm}\|\F(m)-\F_\theta(m)\|_\Y<\epsilon,\qquad
 \sup_{m\in\Adm}\|(D\F(m)-D\F_\theta(m))T\|_{\Lop(\M,\Y)}<\epsilon.
\]
\end{corollary}
\begin{proof}
The set $H_T=\overline{T B_1^\M}$ is compact in $\M$. Apply
Theorem~\ref{thm:compact-directions} with direction set $H_T$ and note that
$\|(D\F(m)-D\F_\theta(m))T\|
 \le\sup_{h\in H_T}\|(D\F(m)-D\F_\theta(m))h\|_\Y$.
\end{proof}

Unlike Theorem~\ref{thm:operator-approx}, this result does not require
$D\F(m)$ to be compact, because the weight $T$ supplies the compactness. It is
used with $T=C^{1/2}$ in Section~\ref{sec:preconditioned-landweber}.

Section~\ref{sec:architectures} writes the reduced-basis class as two
neural-operator architectures, so the results above apply to both.

\section{Generalised DeepONet and generalised Fourier architectures}
\label{sec:architectures}

This section writes the reduced-basis class \eqref{eq:rb} as two
neural-operator architectures. In a generalised DeepONet, the encoder extracts
scalar features from the input, the network predicts output coefficients, and
the decoder combines fixed output elements with those coefficients. In a
generalised Fourier operator, the same class is expressed in Fourier
coordinates.

\subsection{Generalised DeepONet}
\label{sec:generalised-deeponet}

\begin{definition}[Generalised DeepONet]
Let $\ell_1,\ldots,\ell_r\in\M^*$ be bounded linear sensors,
$\beta_1,\ldots,\beta_s\in\Y$ fixed output elements, and
$f_\theta:\mathbb R^r\to\mathbb R^s$ a neural network. Set
\begin{equation}
 \F^{\mathrm{GD}}_\theta(m)
 =\sum_{j=1}^s[f_\theta(\ell_1(m),\ldots,\ell_r(m))]_j\beta_j.
 \label{eq:gdon}
\end{equation}
\end{definition}

The generalised DeepONet formula is another way of writing the operators
\eqref{eq:rb}, with the same derivative when $f_\theta$ is continuously
differentiable.

\begin{proposition}[Exact DeepONet representation]
\label{prop:gdon}
The class \eqref{eq:gdon} equals the class \eqref{eq:rb}. If
$f_\theta\in C^1$, the identification also preserves the Fr\'echet derivative:
\begin{equation}
 D\F^{\mathrm{GD}}_\theta(m)h
 =\sum_{j=1}^s[Df_\theta(\Psi m)\Psi h]_j\beta_j
 =\Phi Df_\theta(\Psi m)\Psi h.
 \label{eq:architecture-derivative}
\end{equation}
\end{proposition}
\begin{proof}
For \eqref{eq:rb}, take $\ell_i=e_i^*\Psi$ and $\beta_j=\Phi e_j$.
Expanding the vector $f_\theta(\Psi m)$ in the standard basis yields
\eqref{eq:gdon}. Conversely, the finitely many bounded sensors define a
bounded encoder, and $c\mapsto\sum_jc_j\beta_j$ is a bounded decoder.
The chain rule gives \eqref{eq:architecture-derivative}.
\end{proof}

The generalised formulation builds on several constructions in the
operator-learning literature. In the original DeepONet, a branch network
maps sampled input values to coefficients, while a trunk network represents
the output functions combined with those coefficients \citep{Lu2021}.
In \eqref{eq:gdon}, $f_\theta$ plays the role of the branch network and the
fixed elements $\beta_j$ play the role of the trunk functions. The input
features are specified more generally by bounded linear functionals
$\ell_i$.

Choosing fixed output functions follows the approach of POD-DeepONet,
where the trunk is constructed from modes obtained by proper orthogonal
decomposition (POD) of training outputs \citep{Lu2022POD}. A related
construction is PCA-Net, which uses principal component analysis (PCA)
to obtain finite-dimensional representations of both inputs and outputs,
and trains a neural network to map between their coordinates
\citep{Bhattacharya2021}. These constructions provide concrete choices
of the input measurements and output elements in \eqref{eq:gdon}.

The specific architecture used in our numerical experiments, denoted by
\emph{GenDeepONet}, adopts empirical PCA/POD representations within this
class. Its encoder extracts coordinates in an empirical input subspace,
$f_\theta$ predicts output coefficients, and the decoder reconstructs the
output using empirical output modes fixed before network training. Fixed
means and affine normalisations can be absorbed into network biases and an
additional constant output element, and they do not change the derivative
formula for the linear parts.

\begin{remark}[Simultaneous approximation by standard DeepONets]
\label{rem:trained-trunk}
In the original DeepONet, the output elements are produced by a trained
trunk network and the sensors are point evaluations. The approximation results extend to this setting under two conditions. First, the sensors must
be bounded on $\mathcal M$, which holds for point evaluations on $C(\overline\Omega)$
but not on $L^2(\Omega)$. Second, the trunk must approximate the fixed output
elements of  \eqref{eq:gdon} in $\mathcal Y$, which is possible whenever network
functions of the output variable are dense in $\mathcal Y$, as for $L^2(\Omega)$ or
$C(\overline\Omega)$. 

Under the respective hypotheses of Section~\ref{sec:approx}, the approximation results extend to standard DeepONets with the additional conditions above.

Fix \(\F_\theta=\Phi f_\theta\Psi\) of the form \eqref{eq:rb}, where \(\Psi=(\ell_1,\ldots,\ell_r)\). For \(\mathcal M=C(X)\) with \(X=\overline\Omega\) compact, take finite partitions of unity \((\varphi_{n,i})_i\) subordinate to covers of vanishing diameter, and choose sample points \(x_{n,i}\) in the corresponding cover sets. Equicontinuity yields
\[
S_nm:=\sum_i m(x_{n,i})\varphi_{n,i}\longrightarrow m
\quad\text{uniformly for }m\in \Adm\cup H.
\]
Thus, setting \(E_nm=(m(x_{n,i}))_i\), the encoders
\[
\Psi_n=\Psi S_n=A_nE_n,
\qquad
(A_n)_{ji}=\ell_j(\varphi_{n,i}),
\]
approximate \(\Psi\) uniformly on \(\Adm\cup H\). This covers Theorem~\ref{thm:compact-directions} and Corollary~\ref{cor:weighted-approximation}, taking \(H=\overline{TB_1^{\mathcal M}}\) in the latter.

For Theorem~\ref{thm:operator-approx}, additionally assume that \(\mathcal M\) is a function space with bounded point evaluations separating its elements. The Hahn--Banach theorem and reflexivity give
\[
\overline{\operatorname{span}\{\delta_x:x\in X\}}
^{\|\cdot\|_{\mathcal M^*}}=\mathcal M^*.
\]
Approximating each \(\ell_j\) by finite linear combinations of point evaluations gives \(\Psi_n=A_nE_n\to\Psi\) in operator norm.

In both cases, continuity of \(f_\theta,Df_\theta\) and
\[
D\F_{\theta,n}(m)=\Phi Df_\theta(\Psi_nm)\Psi_n,
\qquad \F_{\theta,n}=\Phi f_\theta\Psi_n,
\]
give the corresponding value and derivative approximation; \(A_n\) is absorbed into the first branch layer. Finally, replacing \(\beta_j=\Phi e_j\) by trunk functions satisfying \(\|\widetilde\beta_j-\beta_j\|_{\mathcal Y}\le\eta\) changes values and directional derivatives by at most
\[
\eta\|f_\theta(\Psi_nm)\|_1,
\qquad
\eta\|Df_\theta(\Psi_nm)\Psi_nh\|_1,
\]
respectively. These factors are uniformly bounded for \(m\in\Adm\) and \(h\) in the prescribed compact direction set, or in \(B_1^{\mathcal M}\) for Theorem~\ref{thm:operator-approx}. Choosing the sensor and trunk errors successively small completes the approximation.
\end{remark}

The input measurements must be bounded, so that small changes of the input
produce small changes in the encoded coordinates. For $L^2(\Omega)$ inputs, suitable
measurements include inner products with fixed $L^2(\Omega)$ functions,
such as the coordinates used by the PCA encoder. Point samples can also
be used when the input space controls pointwise values, as in
$C(\overline\Omega)$ equipped with the supremum norm.

\subsection{Generalised fully coupled Fourier operators}
\label{sec:generalised-fourier}

Let $V_r,V_s$ be finite-dimensional real Fourier spaces. We use a real
trigonometric basis for each space, consisting of sine and cosine functions
and the constant function when that mode is retained. The coordinate map
$\mathscr C_j:V_j\to\mathbb R^j$ extracts the real coefficients of an
element of $V_j$, while $\mathscr S_j:\mathbb R^j\to V_j$ reconstructs
that element from its coefficients, so that $\mathscr C_j\mathscr S_j=I$ and
$\mathscr S_j\mathscr C_j=I_{V_j}$. This allows the real-valued networks of
Section~\ref{sec:approx} to act on Fourier coefficients.

\begin{definition}[Generalised fully coupled Fourier operator]
Given bounded linear maps $P:\M\to V_r$, $Q:V_s\to\Y$, define
\begin{equation}
 \F^{\mathrm{GF}}_\theta
 =Q\mathscr S_s f_\theta\mathscr C_rP.
 \label{eq:gfno}
\end{equation}
Here $f_\theta:\mathbb R^r\to\mathbb R^s$ is the finite-dimensional
neural network introduced in Section~\ref{sec:approx}. It maps the
Fourier-coordinate vector $\mathscr C_rPm$ to output coefficients, which
$\mathscr S_s$ uses to reconstruct an element of $V_s$. We allow fully
connected layers in $f_\theta$, so each output coefficient may depend on
all retained input coordinates.
\end{definition}
We note that we avoid the name ``generalised FNO'', which has been used for a different
Fourier-based architecture \citep{Ke2026GFNO}.

\begin{proposition}[Exact generalised Fourier representation]
\label{prop:gfno}
If the bounded outer maps $P,Q$ may be chosen freely, the class
\eqref{eq:gfno} equals \eqref{eq:rb}. The operators and their derivatives
agree exactly.
\end{proposition}
\begin{proof}
For $\Phi f_\theta\Psi$, choose
$P=\mathscr S_r\Psi$ and $Q=\Phi\mathscr C_s$. Then
\[
 Q\mathscr S_s f_\theta\mathscr C_rP
 =\Phi(\mathscr C_s\mathscr S_s)f_\theta
       (\mathscr C_r\mathscr S_r)\Psi
 =\Phi f_\theta\Psi.
\]
Conversely, set $\Psi=\mathscr C_rP$ and $\Phi=Q\mathscr S_s$.
These are bounded linear maps, so the chain rule again gives
\eqref{eq:architecture-derivative} with these $\Psi,\Phi$.
\end{proof}

The freedom in $P,Q$ is essential here, since the Fourier maps then act as
coordinate changes around an arbitrary reduced-basis operator. The statement
does not identify every such operator with a model whose encoder is a fixed
Fourier truncation and whose decoder is a fixed Fourier synthesis.

The specific architecture used in our numerical experiments, denoted by
\emph{GenFNO}, uses fixed Fourier truncation and synthesis maps. Its
finite-dimensional network $f_\theta$ is a residual multilayer perceptron
that maps the retained input Fourier coefficients to output Fourier
coefficients. Learning a map between spectral coefficients follows the
Fourier spectral neural operator approach, or SNO(F), of \citet{Fanaskov2023}, and the explicit mode mixing of the higher-order FNO of
\citet{Colagrande2026} is a related construction.

\begin{corollary}[Transfer of simultaneous approximation]
\label{cor:architecture-transfer}
The conclusions of Theorems~\ref{thm:compact-directions} and
\ref{thm:operator-approx}, and of Corollary~\ref{cor:weighted-approximation},
hold with generalised DeepONets \eqref{eq:gdon} and with the broad
generalised Fourier class \eqref{eq:gfno}, under their respective hypotheses.
If $\M=H^s(\mathbb T^d)$ and $\Y=H^{s'}(\mathbb T^d)$, the approximants in these
results can be chosen in \eqref{eq:gfno} with $P$ the orthogonal Fourier
truncation onto $V_r$ and $Q$ the inclusion of $V_s$ into $\Y$.
\end{corollary}
\begin{proof}
Let $\F_\theta=\Phi f_\theta\Psi$ be an approximant from
Theorem~\ref{thm:compact-directions}, Theorem~\ref{thm:operator-approx} or
Corollary~\ref{cor:weighted-approximation}. By Proposition~\ref{prop:gdon} it
equals the generalised DeepONet with $\ell_i=e_i^*\Psi$ and
$\beta_j=\Phi e_j$. By Proposition~\ref{prop:gfno} it equals the element of
\eqref{eq:gfno} with $P=\mathscr S_r\Psi$ and $Q=\Phi\mathscr C_s$. In both
cases the derivatives coincide by \eqref{eq:architecture-derivative}, so the
errors are unchanged.

For the periodic case, let $\{\varphi_k\}$ be the real trigonometric basis,
normalised to be orthonormal in $H^s(\mathbb T^d)$, and let $\Pi_n$ be the
orthogonal projection onto $\operatorname{span}\{\varphi_k:|k|_\infty\le n\}$;
define $\Pi_N$ similarly in $\Y$. Ordering the frequencies by $|k|_\infty$
makes $(\Pi_n)$ a subsequence of the partial-sum projections of an
orthonormal basis. Hence $\|\Pi_n\|=1$, $\Pi_n\to I$ strongly and uniformly on
compact sets, and, since $\Pi_n^*=\Pi_n$, also $\Pi_n^*\to I$ strongly. These
are the only properties of the projections used in
Appendices~\ref{proof:compact-directions} and~\ref{proof:operator-approx}.
Those proofs therefore apply with $P=\Pi_n$ and $Q=\Pi_N$ in the first, and
with $P_n=\Pi_n$ and $Q_N=\Pi_N$ in the second. The encoder $\Psi$ is then the
vector of retained Fourier coefficients and the decoder $\Phi$ is Fourier
synthesis, which is \eqref{eq:gfno} with $P=\Pi_n:\M\to V_r$ and $Q$ the
inclusion $V_s\hookrightarrow\Y$.
\end{proof}

The corollary allows truncations and network sizes to grow and says nothing
about the particular mode counts used in the experiments. The reaction benchmark of Section~\ref{sec:reaction-benchmark} uses finite-grid Fourier coordinates on a bounded square, and this
periodic representation is a computational choice that does not satisfy the
periodic continuum hypotheses above.

\subsection{Comparison with standard FNOs}
\label{sec:fno-comparison}

Standard FNOs, which also appear in our numerical comparisons, arrange their
layers differently from the classes above.

In the Fourier coefficient architecture of
Section~\ref{sec:generalised-fourier}, the finite network $f_\theta$ maps
retained input Fourier coordinates to output coefficients, from which
the decoder reconstructs the output. Its fully connected layers allow
each predicted coefficient to depend on all retained input coordinates.
A standard FNO instead applies successive layers to intermediate
functions, combining Fourier multipliers with nonlinear activations in
physical space \citep{Li2021}. Each multiplier acts separately at each
retained frequency, while the physical-space nonlinearities allow
interactions between frequencies. Appendix~\ref{app:fno-calculus} gives
the standard-FNO layer formulas and their derivative and norm calculations.

\begin{remark}[Frequency coupling in standard FNO layers]
\label{rem:fno-frequency}
Write a standard Fourier layer as
$v_{t+1}=\sigma\bigl(W_tv_t+\mathcal F^{-1}(R_t\cdot\mathcal Fv_t)+b_t\bigr)$,
where $W_t$ acts pointwise on the channels and $R_t(k)$ is the multiplier at
the retained frequency $k$. In Fourier coordinates, the argument of $\sigma$
is $(W_t+R_t(k))\widehat v_t(k)+\widehat b_t(k)$ at each retained frequency,
and $W_t\widehat v_t(k)+\widehat b_t(k)$ at the others. The linear part is therefore block diagonal in $k$. It couples channels but
not frequencies, and coupling between frequencies arises only through the pointwise nonlinearity
in physical space. The generalised Fourier class \eqref{eq:gfno} replaces
this block-diagonal structure by a dense map acting on all retained
coefficients at once, which is what allows it to represent every operator of
the form \eqref{eq:rb}.
\end{remark}

Approximation theory for standard FNOs includes universality for continuous
operators \citep{KovachkiLanthaler2021} and simultaneous approximation of
values and derivatives \citep[Theorem~3.2]{Yao2025}. The latter result
concerns FNOs between periodic Sobolev spaces, gives uniform approximation on
compact input sets, and assumes that the derivative depends continuously on
the input in a Hilbert--Schmidt norm. The space of derivative directions may
be the original input space or a smoother Sobolev space. On the same pair of
spaces, a Hilbert--Schmidt bound also controls the operator norm, but when the
direction space is smoother it controls only directions measured in the
stronger norm, not the full unit ball of the original input space.

In their terminology, a standard FNO trained with derivative information is
a derivative-informed Fourier neural operator (DIFNO). The distinction
concerns the training, not the architecture.

Theorem~\ref{thm:operator-approx}, by contrast, controls the derivative
error of the reduced-basis class over all unit directions of the original
input space, without Hilbert--Schmidt derivatives, and in Hilbert spaces this
also controls the adjoint error in the Landweber update.
\citet[Theorem~3.3]{Yao2025} use derivative approximation to analyse solutions
of strongly convex regularised optimisation problems, whereas we use it for
surrogate Landweber iteration and its stopping rule. The two sets of results
are complementary, and we make no priority claim for simultaneous
approximation of an operator and its derivative.

\section{Preconditioned surrogate Landweber regularisation}
\label{sec:landweber}

In this section, we study preconditioned surrogate Landweber regularisation (PSLR) for the problem 
\begin{equation}\label{eq:model}
	\F(m) = y\,.
\end{equation}
Based on the available data $y^\delta$ satisfying \eqref{eq:observations}, this corresponds to the iteration
\begin{equation}
	m_{k+1}=m_k-\omega\,C\,D\F_\theta(m_k)^*\bigl(\F_\theta(m_k)-y^\delta\bigr),
	\qquad k=0,1,2,\dots,
	\label{eq:preconditioned-surrogate-lw}
\end{equation}
with initial guess $m_0$, step size $\omega>0$ and a self-adjoint, positive
preconditioning operator $C$ on the parameter space $\M$. The employed operator
$\F_\theta$ is a surrogate for the original model or forward operator $\F$ in \eqref{eq:model},
e.g.\ a generalised DeepONet or a member of the generalised Fourier class as discussed before.

Preconditioned Landweber regularisation with the exact
forward map and preconditioners given by negative powers of a Hilbert-scale
generator was studied, e.g., by \citet{Neubauer2000} and
\citet{EggerNeubauer2005}.

Employing a preconditioner $C$ in \eqref{eq:preconditioned-surrogate-lw} offers several advantages. 
Most importantly, in practice the training of the surrogate is based on finitely many directional derivatives $D\F(m)h$ only, and the directions $h$ are typically
sampled from a distribution of smooth functions. A trained network therefore cannot be expected to reproduce $D\F(m)$ in every direction of the
parameter space: in rough directions, which are essentially absent from the
training data, the derivative of the surrogate may be inaccurate however small
the training loss is. If the directions are sampled from a centred Gaussian
measure with trace-class covariance $C'$, the expected squared derivative error
is

\begin{equation}
	\begin{split}
		\mathbb E_{h\sim\mathcal N(0,C')}\|(D\F_\theta(m)-D\F(m))h\|_\Y^2
		&=\|(D\F_\theta(m)-D\F(m))(C')^{1/2}\|_{\mathrm{HS}}^2\\
		&\ge\|(D\F_\theta(m)-D\F(m))(C')^{1/2}\|_{\Lop(\M,\Y)}^2 .
	\end{split}
	\label{eq:training-bound}
\end{equation}
The training thus controls the derivative error weighted by $(C')^{1/2}$, i.e.\ the
derivative in the directions that are smooth with respect to $C'$, but not the
unweighted error $\|D\F_\theta(m)-D\F(m)\|_{\Lop(\M,\Y)}$. A training loss averaged
over parameters $m$ controls this weighted error only on average, whereas the
analysis below requires a bound that is uniform in $m$.

Assume now that the preconditioner $C$ in \eqref{eq:preconditioned-surrogate-lw}
is at least as smoothing as $C'$, in the sense that
\begin{equation}
	\ip{Ch}{h}_\M\le\lambda\ip{C'h}{h}_\M\qquad\text{for all }h\in\M
	\label{eq:smoothing-order}
\end{equation}
for some $\lambda>0$. Applying \eqref{eq:smoothing-order} to
$h=(D\F_\theta(m)-D\F(m))^*y$, $y\in\Y$, yields
\[
\|(D\F_\theta(m)-D\F(m))C^{1/2}\|_{\Lop(\M,\Y)}
\le\lambda^{1/2}\|(D\F_\theta(m)-D\F(m))(C')^{1/2}\|_{\Lop(\M,\Y)} .
\]
The training therefore also controls the derivative error weighted by $C^{1/2}$,
which is the only form in which the derivative error enters the analysis below.
In practice, however, the directions are often taken from a fixed reduced
basis $h_1,\ldots,h_r$, as in Section~\ref{sec:directional-training} and in
our experiments. Then $C'=\sum_{j=1}^r h_j\otimes h_j$ has finite rank, and
\eqref{eq:smoothing-order} cannot hold for an injective $C$: for every $h$
orthogonal to $h_1,\ldots,h_r$ it forces $\ip{Ch}{h}_\M=0$, that is, $Ch=0$.
Training then controls the derivative error only on
$\operatorname{span}\{h_1,\ldots,h_r\}$, while a smoothing $C$ attenuates,
but does not bound, the error in the remaining directions; see
\eqref{eq:supervision-preconditioner-bound} in
Section~\ref{sec:weighted-derivative-error}. The analysis below therefore
assumes a bound on the weighted derivative error rather than deriving it from
the training loss. At the same time, every increment of
\eqref{eq:preconditioned-surrogate-lw} lies in $\ran C$, so a smoothing $C$
keeps the iteration away from the rough directions in which the learned
derivative is least reliable.

\bigskip

After introducing the PSLR \eqref{eq:preconditioned-surrogate-lw}, our aim
is now to show that it is in fact a regularisation method in the classical sense:
\begin{enumerate}
	\item It is well-defined, i.e.\ the iterates $m_k$ are computable.
	\item It is stable with respect to the data $y^\delta$.
	\item It converges towards an exact solution when the noise level $\delta$ and the approximation error of $\F_\theta$ for the true model $\F$ tend to $0$.
\end{enumerate}

Note that stability is a standard result for (exact) Landweber regularisation, see, e.g., \citet{HankeNeubauerScherzer1995, KNS2008}. There it is shown that, under continuity assumptions on $\F$, for each fixed $k$ the iterate $m_k$ depends continuously on $y^\delta$. This result directly carries over to \eqref{eq:preconditioned-surrogate-lw} provided that $\F_\theta$ and the adjoint of its derivative are continuous. 

Of further interest are results on convergence rates, which require
source conditions on $m^\dagger-m_0$; for Hilbert-scale preconditioners and the exact
forward map, their dependence on $C$ is analysed by \citet{Neubauer2000} and
\citet{EggerNeubauer2005}. However, such results for the surrogate case are beyond the scope of the manuscript at hand.

\subsection{Standing assumptions and linearisation error}
\label{sec:preconditioned-landweber}

In this subsection, we make the setting precise and bound the linearisation error of the surrogate. Let $\M$ and $\Y$ be the
separable Hilbert spaces of Section~\ref{sec:framework}. All adjoints are taken with respect to $\ip{\cdot}{\cdot}_\M$.
We group the hypotheses into assumptions on the forward problem and the
surrogate, on the preconditioner and the initial guess, on the nonlinearity, and
on the parameters of the method, and discuss each group directly after it.

\begin{assumption}[Forward problem and surrogate]
	\label{ass:forward}
	The forward operator and its surrogate satisfy $\F,\F_\theta\in C^1(U,\Y)$ for an
	open set $U\subseteq\M$, there exist exact data
	$y^\dagger\in\Y$ and an exact solution $m^\dagger\in U$ of \eqref{eq:model},
	i.e.\ $\F(m^\dagger)=y^\dagger$, and the available data $y^\delta\in\Y$
	satisfy $\|y^\delta-y^\dagger\|_\Y\le\delta$ for a noise level $\delta\ge0$.
\end{assumption}

\begin{remark}
	\label{rem:forward}
	Differentiability is a classical assumption to ensure well-definedness of Landweber iteration; in \eqref{eq:preconditioned-surrogate-lw} it is required for the surrogate $\F_\theta$, since every step evaluates $D\F_\theta(m_k)^*$. Even though we use a surrogate in \eqref{eq:preconditioned-surrogate-lw}, it is reasonable to assume differentiability of $\F$ itself to ensure that $\F$ and its derivative can be approximated by a differentiable surrogate.
\end{remark}

For a self-adjoint, positive and injective operator $C:\M\to\M$, we denote by
\begin{equation}
	\M_C := \ran\left(C^{\frac12}\right),\qquad
	\ip{h}{v}_{\M_C}:=\ip{C^{-1/2}h}{C^{-1/2}v}_\M\,,
	\label{eq:MC-space}
\end{equation}
the space of directions that are smooth with respect to $C$, equipped with the
induced inner product.

\begin{assumption}[Preconditioner and initial guess]
	\label{ass:preconditioner}
	The preconditioning operator $C:\M\to\M$ is self-adjoint, positive, injective
	and compact. For some radius $\rho>0$, the closed ball
	\begin{equation}
		B_\rho^C(m^\dagger):=\{m\in m^\dagger+\M_C:\ \|m-m^\dagger\|_{\M_C}\le\rho\}
		\label{eq:admissible-set}
	\end{equation}
	around the exact solution with respect to the norm of $\M_C$ obeys
	$B_\rho^C(m^\dagger)\subseteq U$, and the initial guess satisfies $m_0\in B_\rho^C(m^\dagger)$, i.e.\
	$m^\dagger-m_0\in\M_C$ and $\|m_0-m^\dagger\|_{\M_C}\le\rho$.
\end{assumption}

\begin{remark}
	\label{rem:preconditioner}
	Self-adjointness, positivity and injectivity of $C$ make $\M_C$ in
	\eqref{eq:MC-space} a Hilbert space, in which $CA^*$ is the adjoint of
	$A|_{\M_C}$ for every $A\in\Lop(\M,\Y)$. Hence
	\eqref{eq:preconditioned-surrogate-lw} is the classical Landweber iteration,
	carried out in $\M_C$. Compactness of $C$ implies that $B_\rho^C(m^\dagger)$ in \eqref{eq:admissible-set} is compact in
	$\M$. This will allow us to employ the previous results on uniform approximation of $\F$ and $D\F$ by the surrogate on $B_\rho^C(m^\dagger)$, the latter in the weighted norm $\|(\cdot)\,C^{1/2}\|_{\Lop(\M,\Y)}$ (Corollary~\ref{cor:weighted-approximation}). In view of \eqref{eq:training-bound},
	$C$ should be at least as smoothing as the covariance $C'$ of the training
	directions in the sense of \eqref{eq:smoothing-order}.
	
	The condition $m_0\in B_\rho^C(m^\dagger)$ requires the initial error to be smooth
	with respect to $C$; this is the price of preconditioning. For $\M=L^2(\Omega)$ on a bounded Lipschitz domain $\Omega\subset\mathbb R^d$, a typical choice is the covariance operator $(Cf)(x)=\int_\Omega k(x-y)f(y)\,\mathrm dy$ of a Mat\'ern kernel $k$ with smoothness $\nu_C$. Then $\M_C$ coincides with $H^{\nu_C+d/2}(\Omega)$,
	with equivalent norms and without boundary conditions
	\citep[Chapter~10]{Wendland2005}. Hence a discontinuous initial error is excluded
	for every $\nu_C>0$, whereas an initial error that is a Mat\'ern sample of
	smoothness $\nu>\nu_C+d/2$ lies in $\M_C$ almost surely \citep{Kanagawa2018}. The set $B_\rho^C(m^\dagger)$ contains $m^\dagger$ and
	depends only on $m^\dagger$, $\rho$ and $C$, but not on the initial guess.
\end{remark}

\begin{assumption}[Nonlinearity]
	\label{ass:surrogate}
	For some
	$0\le\eta<\frac12$, either the forward operator satisfies the local tangential
	cone condition
	\begin{equation}
		\|\F(m)-\F(m^\dagger)-D\F(m)(m-m^\dagger)\|_\Y
		\le\eta\|\F(m)-\F(m^\dagger)\|_\Y
		\label{eq:tcc}\tag{TCC}
	\end{equation}
	for all $m\in B_\rho^C(m^\dagger)$, or the surrogate satisfies
	\begin{equation}
		\|\F_\theta(m)-\F_\theta(m^\dagger)-D\F_\theta(m)(m-m^\dagger)\|_\Y
		\le\eta\|\F_\theta(m)-\F_\theta(m^\dagger)\|_\Y
		\label{eq:tcc-surrogate}\tag{TCC$_\theta$}
	\end{equation}
	for all $m\in B_\rho^C(m^\dagger)$.
\end{assumption}

\begin{remark}
	\label{rem:surrogate-assumption}
	The tangential cone condition is a standard
	assumption for nonlinear Landweber iteration \citep[see, e.g.,][]{HankeNeubauerScherzer1995,KNS2008}. 
	Condition \eqref{eq:tcc} is a property of the physical problem, whereas
	\eqref{eq:tcc-surrogate} is a property of the trained network. Note that the latter
	does not follow from \eqref{eq:tcc} and small forward or derivative training losses,
	since \eqref{eq:tcc} is a relative condition whose right-hand side
	$\eta\|\F(m)-\F(m^\dagger)\|_\Y$ tends to zero as $m\to m^\dagger$, whereas absolute
	approximation errors do not. However, if \eqref{eq:tcc-surrogate} holds for surrogates
	$\F_n$ with a common $\eta$, and if $\F_n(m)\to\F(m)$ and
	$D\F_n(m)(m-m^\dagger)\to D\F(m)(m-m^\dagger)$ for all $m\in B_\rho^C(m^\dagger)$ as
	$n\to\infty$, then letting $n\to\infty$ in \eqref{eq:tcc-surrogate} shows that $\F$
	satisfies \eqref{eq:tcc} with the same $\eta$. That is, when the surrogates and their
	derivatives converge, the surrogate cone condition \eqref{eq:tcc-surrogate} is
	not weaker than the physical one \eqref{eq:tcc} in the limit.
\end{remark}

Now, we quantify the approximation errors of the surrogate, which act as additional noise levels, by
\begin{equation}
	\begin{gathered}
		\epsop:=\sup_{m\in B_\rho^C(m^\dagger)}\|\F_\theta(m)-\F(m)\|_\Y,\\
		\varepsilon_{D,C}:=\sup_{m\in B_\rho^C(m^\dagger)}\|(D\F_\theta(m)-D\F(m))C^{1/2}\|_{\Lop(\M,\Y)} .
	\end{gathered}
	\label{eq:accuracy}
\end{equation}
Here $\epsop$ is the forward error of the surrogate on $B_\rho^C(m^\dagger)$, and
$\varepsilon_{D,C}$ is its derivative error weighted by $C^{1/2}$. Both suprema are finite, since $\F,\F_\theta\in C^1(U,\Y)$ and $B_\rho^C(m^\dagger)$ is compact by Remark~\ref{rem:preconditioner}. Under Assumptions~\ref{ass:forward} and~\ref{ass:preconditioner}, $\F\in C^1(U,\Y)$ and $B_\rho^C(m^\dagger)$ is compact, and hence Corollary~\ref{cor:weighted-approximation} with
$T=C^{1/2}$ provides, for every $\epsilon>0$, a surrogate of the form
\eqref{eq:rb} with $\epsop<\epsilon$ and $\varepsilon_{D,C}<\epsilon$. According to Corollary~\ref{cor:architecture-transfer}, this can also be realized as a generalised
DeepONet or as a member of the generalised Fourier class.

The derivative error also controls the size of the weighted surrogate
derivative. With
$L_C:=\sup_{m\in B_\rho^C(m^\dagger)}\|D\F(m)C^{1/2}\|_{\Lop(\M,\Y)}<\infty$, the triangle
inequality gives
\begin{equation}
	\sup_{m\in B_\rho^C(m^\dagger)}\|D\F_\theta(m)C^{1/2}\|_{\Lop(\M,\Y)}\le L_C+\varepsilon_{D,C}.
	\label{eq:derivative-bound}
\end{equation}
This will be useful for the choice of the step size $\omega$ in \eqref{eq:preconditioned-surrogate-lw}:
\begin{assumption}[Parameters of the method]
	\label{ass:parameters}
	Choose a constant $C_0 > 0$ such that $\sup_{m\in B_\rho^C(m^\dagger)}\|D\F_\theta(m)C^{1/2}\|_{\Lop(\M,\Y)}\le C_0$ and choose the step size $\omega$ such that 
	\begin{equation}
		0<\omega\le C_0^{-2}\,.
		\label{eq:stepsize}
	\end{equation}
\end{assumption}

\begin{remark}
	\label{rem:parameters}
	Under Assumptions~\ref{ass:forward} and~\ref{ass:preconditioner}, the supremum in Assumption~\ref{ass:parameters} is finite, so an admissible constant $C_0$ always exists.
	Any upper bound may be used; by \eqref{eq:derivative-bound}, $C_0=1+L_C$ is
	admissible for every surrogate with $\varepsilon_{D,C}\le1$, independently of
	the particular surrogate.
\end{remark}

\begin{lemma}[Linearisation error of the surrogate]
	\label{lem:linearization}
	Let Assumptions~\ref{ass:forward}, \ref{ass:preconditioner} and
	\ref{ass:surrogate} hold, and let
	\begin{equation}
		\gamma:=
		\begin{cases}
			\rho\,\varepsilon_{D,C} & \text{under \eqref{eq:tcc}},\\[0.5ex]
			0 & \text{under \eqref{eq:tcc-surrogate}},
		\end{cases}
		\label{eq:gamma}
	\end{equation}
	where $\gamma=0$ may be used if both conditions hold. Then,
	for every $m\in B_\rho^C(m^\dagger)$,
	\begin{equation}
		\|D\F_\theta(m)(m-m^\dagger)-(\F_\theta(m)-y^\delta)\|_\Y
		\le\eta\|\F_\theta(m)-y^\delta\|_\Y+(1+\eta)(\delta+\epsop)+\gamma .
		\label{eq:remainder-bound}
	\end{equation}
\end{lemma}

\begin{proof}
	Let $m\in B_\rho^C(m^\dagger)$ and $e:=m-m^\dagger$, so that $e\in\M_C$ and
	$\|e\|_{\M_C}\le\rho$. Moreover,
	\[
	\|y^\delta-\F_\theta(m^\dagger)\|_\Y\le\|y^\delta-y^\dagger\|_\Y
	+\|\F(m^\dagger)-\F_\theta(m^\dagger)\|_\Y\le\delta+\epsop .
	\]
	If the surrogate satisfies \eqref{eq:tcc-surrogate}, define the nonlinear
	remainder $r$ of the surrogate by
	\[
	D\F_\theta(m)e=\F_\theta(m)-\F_\theta(m^\dagger)+r .
	\]
	Then directly $\|r\|_\Y\le\eta\|\F_\theta(m)-\F_\theta(m^\dagger)\|_\Y$ by
	\eqref{eq:tcc-surrogate}, and hence
	\begin{align*}
		\|D\F_\theta(m)e-(\F_\theta(m)-y^\delta)\|_\Y
		&=\|r+(y^\delta-\F_\theta(m^\dagger))\|_\Y\\
		&\le\eta\|\F_\theta(m)-\F_\theta(m^\dagger)\|_\Y+\|y^\delta-\F_\theta(m^\dagger)\|_\Y\\
		&\le\eta\bigl(\|\F_\theta(m)-y^\delta\|_\Y
		+\|y^\delta-\F_\theta(m^\dagger)\|_\Y\bigr)+\|y^\delta-\F_\theta(m^\dagger)\|_\Y\\
		&\le\eta\|\F_\theta(m)-y^\delta\|_\Y+(1+\eta)(\delta+\epsop) .
	\end{align*}
	If the forward operator satisfies \eqref{eq:tcc}, substituting
	$\F_\theta=\F+(\F_\theta-\F)$, $D\F_\theta=D\F+(D\F_\theta-D\F)$ and
	$y^\delta=\F(m^\dagger)+(y^\delta-y^\dagger)$ and using the triangle inequality
	gives
	\begin{align*}
		&\|D\F_\theta(m)e-(\F_\theta(m)-y^\delta)\|_\Y\\
		&\quad=\Bigl\|(D\F_\theta(m)-D\F(m))e
		+D\F(m)e-\bigl(\F(m)-\F(m^\dagger)\bigr)\\
		&\qquad\quad{}+\bigl(\F(m)-\F_\theta(m)\bigr)
		+\bigl(y^\delta-y^\dagger\bigr)\Bigr\|_\Y\\
		&\quad\le\|(D\F_\theta(m)-D\F(m))e\|_\Y
		+\bigl\|D\F(m)e-\bigl(\F(m)-\F(m^\dagger)\bigr)\bigr\|_\Y\\
		&\qquad\quad{}+\|\F(m)-\F_\theta(m)\|_\Y
		+\|y^\delta-y^\dagger\|_\Y\\
		&\quad\le\rho\,\varepsilon_{D,C}
		+\eta\|\F(m)-\F(m^\dagger)\|_\Y+\epsop+\delta\\
		&\quad\le\rho\,\varepsilon_{D,C}
		+\eta\bigl(\|\F(m)-\F_\theta(m)\|_\Y
		+\|\F_\theta(m)-y^\delta\|_\Y
		+\|y^\delta-y^\dagger\|_\Y\bigr)+\epsop+\delta\\
		&\quad\le\eta\|\F_\theta(m)-y^\delta\|_\Y
		+(1+\eta)(\delta+\epsop)+\rho\,\varepsilon_{D,C} .
	\end{align*}
	Here the first term was bounded by the weighted derivative error,
	\[
	\|(D\F_\theta(m)-D\F(m))e\|_\Y
	=\|(D\F_\theta(m)-D\F(m))C^{1/2}C^{-1/2}e\|_\Y
	\le\varepsilon_{D,C}\|e\|_{\M_C}\le\rho\,\varepsilon_{D,C},
	\]
	the second by
	\eqref{eq:tcc}, and the last two by the definition of $\epsop$ and by
	Assumption~\ref{ass:forward}.
\end{proof}

Lemma~\ref{lem:linearization} is the only consequence of \eqref{eq:tcc} or
\eqref{eq:tcc-surrogate} used in the sequel, so that both cases are treated
simultaneously from now on. 

\subsection{Stopping and monotonicity}
\label{sec:stopping}

Let us now derive an appropriate discrepancy principle for stopping the iteration \eqref{eq:preconditioned-surrogate-lw}. Afterwards, we will prove that all iterates up to the corresponding stopping index lie in $B_\rho^C(m^\dagger)$ and are hence well-defined.

As a rationale for the stopping criterion, note that 
\[
\|\F_\theta(m^\dagger)-y^\delta\|_\Y
\le\|\F_\theta(m^\dagger)-\F(m^\dagger)\|_\Y+\|y^\dagger-y^\delta\|_\Y
\le\epsop+\delta\,.
\]
From this point of view, the effective noise level in the surrogate equation $\F_\theta(m)=y^\delta$ is $\epsop+\delta$. However, the approximation of $D\F$ may cause further errors, as seen in Lemma~\ref{lem:linearization}. Therefore, we define
\begin{equation}\label{eq:effective_noiselevel}
	\Delta:=\delta+\epsop+\frac{\gamma}{1+\eta}
\end{equation}
with $\gamma$ as in Lemma~\ref{lem:linearization}, and stop the iteration using the discrepancy principle	
\begin{subequations}\label{eqs:discrepancy_principle}
	\begin{equation}
		k_*:=\inf\{k\ge0:\ \|\F_\theta(m_k)-y^\delta\|_\Y\le\tau\Delta\}
		\label{eq:stopping}
	\end{equation}
	with the convention $\inf \emptyset = \infty$ and a parameter $\tau > 0$ satisfying
	\begin{equation}\label{eq:tau}
		\tau>\frac{2(1+\eta)}{1-2\eta}\,.
	\end{equation}
\end{subequations}
The lower bound on $\tau$ in \eqref{eq:tau} is the usual condition of the
discrepancy principle for nonlinear (exact) Landweber iteration \citep[see, e.g.,][]{HankeNeubauerScherzer1995}.
For the discrepancy principle to work, we have to assume that $\Delta>0$. Note that $\Delta=0$ corresponds to exact data and a surrogate that coincides with $\F$ on $B_\rho^C(m^\dagger)$, i.e.\ \eqref{eq:preconditioned-surrogate-lw} reduces to the preconditioned Landweber iteration with exact forward operator and exact data, which is not stopped.
Moreover, the approximation error $\varepsilon_{D,C}$ of the derivative enters the
discrepancy principle only through $\gamma$, i.e.\ only if \eqref{eq:tcc} is used;
under \eqref{eq:tcc-surrogate}, $\Delta=\delta+\epsop$.

In the following, we denote the error of the $k$-th iterate by $e_k:=m_k-m^\dagger$.
Following the standard analysis of Landweber iteration, we first show that a
single step from an iterate in $B_\rho^C(m^\dagger)$ decreases the error in the
norm of $\M_C$ as long as the discrepancy principle has not stopped the
iteration. Since Lemma~\ref{lem:linearization} and the step-size condition
\eqref{eq:stepsize} are only available on $B_\rho^C(m^\dagger)$, we then show by
induction that the iteration does not leave $B_\rho^C(m^\dagger)$ before it is
stopped.

\begin{lemma}[Recursive error bound, monotonicity]
	\label{lem:descent}
	Let Assumptions~\ref{ass:forward}--\ref{ass:parameters} hold, let $\Delta$ be
	given by \eqref{eq:effective_noiselevel} and let $\tau$ satisfy \eqref{eq:tau}.
	If $m_k\in B_\rho^C(m^\dagger)$ and $\|\F_\theta(m_k)-y^\delta\|_\Y>\tau\Delta$,
	then $m_{k+1}$ from \eqref{eq:preconditioned-surrogate-lw} satisfies
	$e_{k+1}\in\M_C$ and
	\begin{equation}
		\|e_{k+1}\|_{\M_C}^2-\|e_k\|_{\M_C}^2
		\le-c_\tau\,\omega\,\|\F_\theta(m_k)-y^\delta\|_\Y^2<0
		\label{eq:descent}
	\end{equation}
	with $c_\tau:=1-2\eta-2(1+\eta)/\tau>0$.
\end{lemma}

\begin{proof}
	Since $m_k\in B_\rho^C(m^\dagger)\subseteq U$, the iterate $m_{k+1}$ is defined,
	and $e_k\in\M_C$. Moreover, $Cz=C^{1/2}(C^{1/2}z)\in\ran(C^{1/2})=\M_C$ for
	every $z\in\M$, so that also
	\[
	e_{k+1}=e_k-\omega CD\F_\theta(m_k)^*(\F_\theta(m_k)-y^\delta)\in\M_C .
	\]
	Expanding the squared error in the norm of $\M_C$ gives
	\begin{equation}
		\begin{split}
			\|e_{k+1}\|_{\M_C}^2
			={}&\|e_k\|_{\M_C}^2
			-2\omega\ip{e_k}{CD\F_\theta(m_k)^*(\F_\theta(m_k)-y^\delta)}_{\M_C}\\
			&+\omega^2\|CD\F_\theta(m_k)^*(\F_\theta(m_k)-y^\delta)\|_{\M_C}^2 .
		\end{split}
		\label{eq_monotone}
	\end{equation}
	As noted in Remark~\ref{rem:preconditioner}, $CD\F_\theta(m_k)^*$ is the adjoint
	of $D\F_\theta(m_k)|_{\M_C}$ in $\M_C$. For the norm, we use $C^{-1/2}C=C^{1/2}$
	and, since $C^{1/2}$ is self-adjoint,
	\[
	\|C^{1/2}D\F_\theta(m_k)^*\|_{\Lop(\Y,\M)}
	=\|(D\F_\theta(m_k)C^{1/2})^*\|_{\Lop(\Y,\M)}
	=\|D\F_\theta(m_k)C^{1/2}\|_{\Lop(\M,\Y)} .
	\]
	Hence
	\begin{align*}
		\ip{e_k}{CD\F_\theta(m_k)^*(\F_\theta(m_k)-y^\delta)}_{\M_C}
		&=\ip{D\F_\theta(m_k)e_k}{\F_\theta(m_k)-y^\delta}_\Y,\\
		\|CD\F_\theta(m_k)^*(\F_\theta(m_k)-y^\delta)\|_{\M_C}
		&=\|C^{1/2}D\F_\theta(m_k)^*(\F_\theta(m_k)-y^\delta)\|_\M\\
		&\le\|D\F_\theta(m_k)C^{1/2}\|_{\Lop(\M,\Y)}\,\|\F_\theta(m_k)-y^\delta\|_\Y\\
		&\le C_0\|\F_\theta(m_k)-y^\delta\|_\Y ,
	\end{align*}
	where the last step uses Assumption~\ref{ass:parameters} and
	$m_k\in B_\rho^C(m^\dagger)$. Together with $\omega C_0^2\le1$ from
	\eqref{eq:stepsize}, we obtain
	\[
	\|e_{k+1}\|_{\M_C}^2-\|e_k\|_{\M_C}^2
	\le-2\omega\ip{D\F_\theta(m_k)e_k}{\F_\theta(m_k)-y^\delta}_\Y+\omega\|\F_\theta(m_k)-y^\delta\|_\Y^2 .
	\]
	Furthermore, the Cauchy--Schwarz inequality and Lemma~\ref{lem:linearization}
	with $m=m_k$ show that the inner-product term satisfies
	\begin{align*}
		&{-}2\omega\ip{D\F_\theta(m_k)e_k}{\F_\theta(m_k)-y^\delta}_\Y\\
		&\quad=-2\omega\|\F_\theta(m_k)-y^\delta\|_\Y^2
		-2\omega\ip{D\F_\theta(m_k)e_k-(\F_\theta(m_k)-y^\delta)}{\F_\theta(m_k)-y^\delta}_\Y\\
		&\quad\le-2\omega\|\F_\theta(m_k)-y^\delta\|_\Y^2
		+2\omega\,\|D\F_\theta(m_k)e_k-(\F_\theta(m_k)-y^\delta)\|_\Y\,\|\F_\theta(m_k)-y^\delta\|_\Y\\
		&\quad\le-2\omega\|\F_\theta(m_k)-y^\delta\|_\Y^2\\
		&\qquad\quad{}+2\omega\bigl(\eta\|\F_\theta(m_k)-y^\delta\|_\Y+(1+\eta)(\delta+\epsop)+\gamma\bigr)\|\F_\theta(m_k)-y^\delta\|_\Y .
	\end{align*}
	Substituting this bound and using $(1+\eta)(\delta+\epsop)+\gamma=(1+\eta)\Delta$
	by \eqref{eq:effective_noiselevel} yields
	\begin{align*}
		&\|e_{k+1}\|_{\M_C}^2-\|e_k\|_{\M_C}^2\\
		&\quad\le-\omega\|\F_\theta(m_k)-y^\delta\|_\Y^2
		+2\omega\bigl(\eta\|\F_\theta(m_k)-y^\delta\|_\Y+(1+\eta)\Delta\bigr)\|\F_\theta(m_k)-y^\delta\|_\Y\\
		&\quad=-\omega\|\F_\theta(m_k)-y^\delta\|_\Y
		\bigl[(1-2\eta)\|\F_\theta(m_k)-y^\delta\|_\Y-2(1+\eta)\Delta\bigr].
	\end{align*}
	By assumption, $\Delta<\frac1\tau\|\F_\theta(m_k)-y^\delta\|_\Y$, and hence
	\begin{align*}
		(1-2\eta)\|\F_\theta(m_k)-y^\delta\|_\Y-2(1+\eta)\Delta
		&>\Bigl(1-2\eta-\frac{2(1+\eta)}{\tau}\Bigr)\|\F_\theta(m_k)-y^\delta\|_\Y\\
		&=c_\tau\|\F_\theta(m_k)-y^\delta\|_\Y ,
	\end{align*}
	where $c_\tau>0$ is equivalent to \eqref{eq:tau}, since $1-2\eta>0$. Finally,
	$\gamma\ge0$ implies $\Delta\ge0$, so that
	$\|\F_\theta(m_k)-y^\delta\|_\Y>\tau\Delta\ge0$, which yields \eqref{eq:descent}
	including the strict inequality.
\end{proof}

\begin{theorem}[Well-definedness and finite stopping]
	\label{thm:stopping}
	Let Assumptions~\ref{ass:forward}--\ref{ass:parameters} hold and let $k_*$ be
	given by the discrepancy principle \eqref{eqs:discrepancy_principle}. Then the
	iterates of \eqref{eq:preconditioned-surrogate-lw} satisfy
	\begin{equation}
		m_k\in B_\rho^C(m^\dagger)\subseteq U\qquad\text{for all integers }0\le k\le k_*,
		\label{eq:well-defined}
	\end{equation}
	i.e.\ the iteration is well defined up to the stopping index, and
	\eqref{eq:descent} holds for all $k<k_*$. If, in addition, $\Delta>0$, then the
	stopping index is finite.  %
\end{theorem}

\begin{proof}
	We argue by induction. By Assumption~\ref{ass:preconditioner},
	$m_0\in B_\rho^C(m^\dagger)$. Let $k<k_*$ with $m_k\in B_\rho^C(m^\dagger)$. By
	\eqref{eq:stopping}, $\|\F_\theta(m_k)-y^\delta\|_\Y>\tau\Delta$, so
	Lemma~\ref{lem:descent} applies and yields $e_{k+1}\in\M_C$ with
	$\|e_{k+1}\|_{\M_C}<\|e_k\|_{\M_C}\le\rho$, that is,
	$m_{k+1}\in B_\rho^C(m^\dagger)$. This proves \eqref{eq:well-defined}, and
	Lemma~\ref{lem:descent} gives \eqref{eq:descent} for every $k<k_*$.
	
	Now let $\Delta>0$ and let $N\ge1$ be an integer with $N\le k_*$. Summing
	\eqref{eq:descent} from $k=0$ to $k=N-1$ gives
	\[
	c_\tau\omega\sum_{k=0}^{N-1}\|\F_\theta(m_k)-y^\delta\|_\Y^2
	\le\|e_0\|_{\M_C}^2-\|e_N\|_{\M_C}^2
	\le\|e_0\|_{\M_C}^2 .
	\]
	Since every residual before stopping satisfies
	$\|\F_\theta(m_k)-y^\delta\|_\Y^2>\tau^2\Delta^2$ by \eqref{eq:stopping},
	\[
	Nc_\tau\omega\tau^2\Delta^2<\|e_0\|_{\M_C}^2 .
	\]
	Hence every integer $1\le N\le k_*$ is bounded by
	$\|e_0\|_{\M_C}^2/(c_\tau\omega\tau^2\Delta^2)$, which proves finiteness of $k_*$.
\end{proof}

\begin{remark}[Variable step sizes]
	\label{rem:step-sizes}
	The step size enters the proof of Lemma~\ref{lem:descent} only through
	\[
	\omega\|C^{1/2}D\F_\theta(m_k)^*(\F_\theta(m_k)-y^\delta)\|_\M^2\le\|\F_\theta(m_k)-y^\delta\|_\Y^2,
	\]
	which is
	\[
	\omega^2\|CD\F_\theta(m_k)^*(\F_\theta(m_k)-y^\delta)\|_{\M_C}^2\le\omega\|\F_\theta(m_k)-y^\delta\|_\Y^2\,,
	\]
	and, for finite
	stopping, through a positive lower bound. Lemma~\ref{lem:descent} and Theorem~\ref{thm:stopping} therefore remain valid for variable step sizes $\omega_k\in[\omega_{\min},C_0^{-2}]$ with some fixed $\omega_{\min}>0$. %
\end{remark}

\subsection{Regularisation convergence as the noise and surrogate errors vanish}
\label{subsec:regularization_convergence}

Theorem~\ref{thm:stopping} shows that, for a fixed noise level and a fixed
surrogate, the iteration is well defined up to the stopping index, the error
decreases monotonically before stopping, and the stopping index is finite. We now
investigate the limiting behaviour of the stopped iterates as both the noise level
and the surrogate errors tend to zero. To this end, we consider sequences of data
and surrogates.

\begin{assumption}[Sequences of data and surrogates]
	\label{ass:sequences}
	For every $n\in\mathbb N$, data $y^{\delta_n}$ and a surrogate
	$\F_n:=\F_{\theta_n}$ are given such that
	Assumptions~\ref{ass:forward}--\ref{ass:parameters} hold for
	$(\F_n,y^{\delta_n})$ with noise level $\delta_n$, step size $\omega_n$ and
	constant $C_0^{(n)}$, where $C$, $m^\dagger$, $\rho$, $m_0$ and $\eta$, and
	hence $B_\rho^C(m^\dagger)$, are the same for all $n$. For each $n$, either \eqref{eq:tcc} or
	\eqref{eq:tcc-surrogate} may be used.
\end{assumption}

\begin{remark}
	\label{rem:sequences}
	Apart from $C$, $m^\dagger$, $\rho$ and $m_0$, only $\eta$ has to be the same
	for all $n$. In particular, the step sizes $\omega_n$ and the constants
	$C_0^{(n)}$ may change with $n$. Under \eqref{eq:tcc}, a common step size can
	nevertheless be chosen for all surrogates with $\varepsilon_{D,C}\le1$, cf.\ Remark~\ref{rem:parameters}.
\end{remark}

For each $n$, we denote by $\epsop^{(n)}$, $\varepsilon_{D,C}^{(n)}$, $\gamma_n$
and $\Delta_n$ the quantities \eqref{eq:accuracy}, \eqref{eq:gamma} and
\eqref{eq:effective_noiselevel} for $\F_n$ and $y^{\delta_n}$, by $m_k^{(n)}$ the
iterates of \eqref{eq:preconditioned-surrogate-lw} with $\F_n$, $y^{\delta_n}$
and $\omega_n$, and by $\widehat m_n:=m_{k_n}^{(n)}$ the iterate stopped by the
discrepancy principle \eqref{eqs:discrepancy_principle} at index $k_n$.

\begin{theorem}[Regularisation property of PSLR]
	\label{thm:convergence}
	\label{thm:surrogate_landweber_regularization}
	Let Assumption~\ref{ass:sequences} hold, let the iterations be stopped by the
	discrepancy principle \eqref{eqs:discrepancy_principle} with a common $\tau$
	satisfying \eqref{eq:tau}, and assume $\Delta_n>0$ for all $n$ as well as
	\[
	\delta_n\to0,\qquad\epsop^{(n)}\to0,\qquad\gamma_n\to0 .
	\]
	Then %
	every subsequence of $(\widehat m_n)$ has a further subsequence
	$(\widehat m_{n_j})$ for which there is a solution
	$\bar m\in B_\rho^C(m^\dagger)$ of \eqref{eq:model} such that
	\[
	\widehat m_{n_j}\to\bar m\quad\text{in }\M
	\]
	If $m^\dagger$ is the only solution of $\F(m)=y^\dagger$ in
	$B_\rho^C(m^\dagger)$, then the whole sequence converges, i.e.\
	$\widehat m_n\to m^\dagger$ in $\M$.%
	Consequently, the PSLR \eqref{eq:preconditioned-surrogate-lw} together with the discrepancy principle \eqref{eqs:discrepancy_principle} is a regularisation method
	with respect to the combined data and surrogate errors.
\end{theorem}

\begin{proof}
	For every $n$, Theorem~\ref{thm:stopping} ensures that $k_n<\infty$ and that
	all iterates up to and including $k_n$ lie in $B_\rho^C(m^\dagger)$; in
	particular, $\widehat m_n\in B_\rho^C(m^\dagger)$. At the stopping index, the
	discrepancy principle \eqref{eq:stopping} gives
	\[
	\|\F_n(\widehat m_n)-y^{\delta_n}\|_\Y\le\tau\Delta_n .
	\]
	The error of the stopped iterate with respect to the exact forward equation
	can be estimated by inserting the surrogate prediction and the noisy datum:
	\begin{align}
		\|\F(\widehat m_n)-y^\dagger\|_\Y
		&\le\|\F(\widehat m_n)-\F_n(\widehat m_n)\|_\Y
		+\|\F_n(\widehat m_n)-y^{\delta_n}\|_\Y
		+\|y^{\delta_n}-y^\dagger\|_\Y\nonumber\\
		&\le\epsop^{(n)}+\tau\Delta_n+\delta_n
		=(1+\tau)\bigl(\delta_n+\epsop^{(n)}\bigr)+\frac{\tau\gamma_n}{1+\eta},\label{eq:true_residual_stopped_iterate}
	\end{align}
	where the first term is bounded by $\epsop^{(n)}$ because
	$\widehat m_n\in B_\rho^C(m^\dagger)$, and the last one by
	Assumption~\ref{ass:forward}. Since $\delta_n\to0$, $\epsop^{(n)}\to0$ and
	$\gamma_n\to0$, estimate \eqref{eq:true_residual_stopped_iterate} shows that
	$\|\F(\widehat m_n)-y^\dagger\|_\Y\to0$. Thus,
	although the discrepancy principle is formulated using the learned surrogate
	$\F_n$ and the noisy data $y^{\delta_n}$, the corresponding stopped iterates
	asymptotically satisfy the exact forward equation.
	
	It remains to relate this vanishing residual to convergence of the unknowns.
	Since every $\widehat m_n$ belongs to $B_\rho^C(m^\dagger)$, the sequence
	$(\widehat m_n-m^\dagger)$ is bounded by $\rho$ in the Hilbert space $\M_C$.
	Consider an arbitrary subsequence of $(\widehat m_n)$. It possesses a further
	subsequence, still denoted for convenience by $(\widehat m_{n_j})$, such that
	$\widehat m_{n_j}-m^\dagger\rightharpoonup\bar m-m^\dagger$ in $\M_C$
	for some $\bar m\in m^\dagger+\M_C$. By weak lower semicontinuity of the norm,
	\[
	\|\bar m-m^\dagger\|_{\M_C}\le\liminf_{j\to\infty}\|\widehat m_{n_j}-m^\dagger\|_{\M_C}\le\rho,
	\]
	that is, $\bar m\in B_\rho^C(m^\dagger)$. Since $B_\rho^C(m^\dagger)$ is compact
	in $\M$ (Remark~\ref{rem:preconditioner}), the embedding
	$\M_C\hookrightarrow\M$ is compact and hence maps weakly convergent sequences
	to strongly convergent ones, so that $\widehat m_{n_j}\to\bar m$ in $\M$.
	The assumption $\F\in C^1(U,\Y)$ implies continuity of $\F$, and therefore
	$\F(\widehat m_{n_j})\to\F(\bar m)$.
	On the other hand, \eqref{eq:true_residual_stopped_iterate} gives
	$\F(\widehat m_{n_j})\to y^\dagger$, i.e.\ $\bar m\in B_\rho^C(m^\dagger)$ is a solution of
	$\F(m)=y^\dagger$.
	
	Finally, suppose that $m^\dagger$ is the only solution of $\F(m)=y^\dagger$
	in $B_\rho^C(m^\dagger)$. Then $\bar m=m^\dagger$ for every subsequence
	considered above, i.e.\ every subsequence of $(\widehat m_n)$ has a further
	subsequence that converges to $m^\dagger$ in $\M$. Since the limit does not depend on the
	subsequence, the whole sequence converges, $\widehat m_n\to m^\dagger$ in
	$\M$.
\end{proof}

\begin{remark}
	Note that convergence $\lim_{n \to \infty} \varepsilon_{D,C}^{(n)} = 0$ of the derivative approximation error $\varepsilon_{D,C}^{(n)}$ from \eqref{eq:accuracy} is required only if the tangential cone condition is assumed for the physical model \eqref{eq:tcc}. If the surrogates $\F_n$ themselves satisfy \eqref{eq:tcc-surrogate}, the derivative approximation errors $\varepsilon_{D,C}^{(n)}$ do not need to converge to $0$.
\end{remark}

\section{Derivative-informed training and parameter-space geometry}
\label{sec:training-preconditioning}

The approximation results of Sections~\ref{sec:approx}
and~\ref{sec:architectures} are existence results and do not say how the
derivative accuracy required in Section~\ref{sec:landweber} can be obtained
from finite training data. In practice, derivative information is supplied
through finitely many directional derivatives. The directions used in
training determine which parts of the derivative are controlled, whereas the
parameter-space geometry used when the surrogate is employed in the inversion determines how the remaining
errors are weighted in the Landweber update. This section makes the
distinction precise and connects it with the choices made in the
experiments.

\subsection{Population and reduced derivative-informed training}
\label{sec:directional-training}

Derivative-informed operator learning augments observations of the forward map
with information about its Fr\'echet derivative.  Following the standard
DINO formulation, let $\mu$ be a probability measure on the parameter space
$\M$ describing the distribution of inputs used for surrogate construction,
and let $\mathcal{H}\hookrightarrow\M$ be a Hilbert space of derivative directions.
A natural population objective is the weighted Sobolev risk
\begin{equation}
    \mathscr R_{\mu,\mathcal{H}}(\theta)
    =
    \mathbb E_{m\sim\mu}
    \left[
        \|\F_\theta(m)-\F(m)\|_\Y^2
        +
        \lambda_J
        \|D\F_\theta(m)-D\F(m)\|_{\mathrm{HS}(\mathcal{H},\Y)}^2
    \right],
    \label{eq:population-dino-risk}
\end{equation}
where $\lambda_J\geq0$ controls the relative contribution of derivative
accuracy, and $\|A\|_{\mathrm{HS}(\mathcal H,\Y)}^2=\sum_{j}\|Ae_j\|_\Y^2$,
for any orthonormal basis $(e_j)$ of $\mathcal H$, is the Hilbert--Schmidt
norm of a linear operator $A:\mathcal H\to\Y$. We write $\|\cdot\|_{\mathrm{HS}}$
for $\|\cdot\|_{\mathrm{HS}(\M,\Y)}$. The risk \eqref{eq:population-dino-risk} is the first-order weighted
Sobolev ($W^{1,2}_\mu$) formulation used in derivative-informed operator
learning
\citep{DINO2024,Yao2025}. In particular, \citet[Theorem~4.3]{Yao2025} establish universal approximation of FNOs in the corresponding $W^{1,2}_{\mu,\delta}$ norm under suitable assumptions, and their Section~5.3 compresses the Hilbert--Schmidt derivative term using reduced bases. The formulation below follows this population-to-reduced-loss
viewpoint and adapts it to the finite set of physical derivative directions
used in our experiments.

The measure $\mu$ determines the region of parameter space over which
approximation accuracy is emphasised. In particular, $\mathcal{H}$ determines the
directions in which derivative accuracy is measured.

A particularly useful choice for covariance-based dimension reduction is
obtained from the covariance of the data-generating measure.  Assuming
$\mu$ has finite second moment, let
\begin{equation}
    \bar m
    =
    \mathbb E_\mu[m],
    \qquad
    \Gamma_\mu v
    =
    \mathbb E_\mu
    \left[
        \langle v,m-\bar m\rangle_\M (m-\bar m)
    \right].
    \label{eq:training-covariance}
\end{equation}
Let $\Gamma_\mu\psi_j=\kappa_j\psi_j$ with $\langle\psi_j,\psi_k\rangle_\M=\delta_{jk}$,
and define the covariance-induced direction space $\M_{\Gamma_\mu} := \operatorname{ran}(\Gamma_\mu^{1/2})$
with its induced Hilbert norm.  The vectors $h_j=\sqrt{\kappa_j}\,\psi_j$
form an orthonormal basis of $\M_{\Gamma_\mu}$ on the support of
$\Gamma_\mu$.  Consequently, writing
\begin{equation}
    B_\theta(m) := D\F_\theta(m)-D\F(m)
    \label{eq:derivative-error-operator}
\end{equation}
for the derivative error, the derivative term in \eqref{eq:population-dino-risk} can be written as
\begin{equation}
    \|B_\theta(m)\|_{\mathrm{HS}(\M_{\Gamma_\mu},\Y)}^2
    =\sum_{j\geq1}\|B_\theta(m)h_j\|_\Y^2
    =\sum_{j\geq1}\kappa_j\|B_\theta(m)\psi_j\|_\Y^2
    =\|B_\theta(m)\Gamma_\mu^{1/2}\|_{\mathrm{HS}(\M,\Y)}^2.
    \label{eq:covariance-weighted-HS}
\end{equation}
When $\mu$ is Gaussian, $\M_{\Gamma_\mu}$ is the corresponding
Cameron--Martin space, although no Gaussian assumption is required for the
covariance-weighted representation above.

Computing the full Hilbert--Schmidt term in
\eqref{eq:covariance-weighted-HS} is generally infeasible after
high-dimensional discretisation.  A standard approximation is therefore to
retain a finite-dimensional subspace.  For the leading $r$ covariance modes,
define
\begin{equation}
    \Gamma_{\mu,r}
    =
    \sum_{j=1}^{r}
    \kappa_j\psi_j\otimes\psi_j .
    \label{eq:truncated-training-covariance}
\end{equation}
The corresponding reduced derivative error is
\begin{equation}
    \|B_\theta(m)\Gamma_{\mu,r}^{1/2}\|_{\mathrm{HS}}^2
    =
    \sum_{j=1}^{r}
    \|B_\theta(m)h_j\|_\Y^2.
    \label{eq:reduced-HS-loss}
\end{equation}
This is the input-reduced Hilbert--Schmidt derivative loss used in
reduced-basis derivative-informed learning.  The discarded part of \eqref{eq:covariance-weighted-HS} lies between $0$
and $\|B_\theta(m)\|_{\Lop(\M,\Y)}^2\sum_{j>r}\kappa_j$, so it is small when
the eigenvalues of $\Gamma_\mu$ decay quickly.

In practice, both the truncated covariance
$\Gamma_{\mu,r}$ in \eqref{eq:truncated-training-covariance} and the
expectation over $\mu$ in \eqref{eq:population-dino-risk} must be estimated
from finite samples.  Let
$m_1,\ldots,m_N$ be independent draws from $\mu$, and let
$h_1,\ldots,h_r$ denote the retained physical derivative directions obtained
from a covariance approximation.  The empirical objective takes the form
\begin{equation}
    \mathcal L_{N,r}(\theta)
    =
    \frac1N\sum_{i=1}^{N}\Bigl[
    \|\F_\theta(m_i)-\F(m_i)\|_\Y^2
    +\frac{\lambda_J}{r}\sum_{j=1}^{r}\|B_\theta(m_i)h_j\|_\Y^2\Bigr].
    \label{eq:empirical-reduced-dino-loss}
\end{equation}
The derivative term averages over the $r$ directions instead of summing,
so that the derivative error in a single direction drawn uniformly at random
is an unbiased estimate of it. Our implementation uses one such direction for
each occurrence of a training input (Section~\ref{sec:reaction-learning}).
Compared with \eqref{eq:population-dino-risk}, the average only rescales
$\lambda_J$ by $1/r$.

With $h_j=\sqrt{\kappa_j}\,\psi_j$, the derivative term is
\eqref{eq:reduced-HS-loss} scaled by $1/r$; the eigenvalues $\kappa_j$ are
carried by the lengths of the $h_j$, since
$\Gamma_{\mu,r}=\sum_{j=1}^r h_j\otimes h_j$. The same holds for any fixed
directions $h_1,\ldots,h_r$, orthogonal or not:
$\sum_{j=1}^r\|B_\theta(m)h_j\|_\Y^2=\|B_\theta(m)\widehat\Gamma_r^{1/2}\|_{\mathrm{HS}}^2$
with $\widehat\Gamma_r=\sum_{j=1}^r h_j\otimes h_j$, which plays the role of
$\Gamma_{\mu,r}$. Finite directional
supervision is therefore a finite-rank approximation of a
covariance-weighted Hilbert--Schmidt derivative loss.

A small value of \eqref{eq:population-dino-risk} gives mean-square control
with respect to $\mu$, and at each parameter its Hilbert--Schmidt term
controls the operator norm on $\mathcal{H}$. Neither the population risk nor
its finite-sample version gives, without further assumptions, the uniform
derivative bound over the admissible set required in
Section~\ref{sec:landweber}.

\subsection{Choosing the preconditioner after training}
\label{sec:parameter-geometry}

Section~\ref{sec:landweber} introduced the preconditioned
iteration \eqref{eq:preconditioned-surrogate-lw} as the setting in which the
approximation and regularisation results combine into
Theorem~\ref{thm:convergence}. Here we consider its practical role. The geometry of the parameter space is
part of the formulation of the inverse method, not a property of the forward
problem, and preconditioned Landweber methods have a long history in inverse
problems \citep{Neubauer2000,EggerNeubauer2005,Egger2005}. Our interest is in
how this choice interacts with the approximation error of a learned forward
operator.

Equivalently, \eqref{eq:preconditioned-surrogate-lw} is obtained by applying
Landweber iteration in $\M$ to $z\longmapsto \F_\theta(m_0+C^{1/2}z)$
and transforming back through $m=m_0+C^{1/2}z$. The operator $C$ is
therefore more than a device to stabilise a learned model, since it specifies
the geometry in which increments are measured. When $C$ is smoothing, $\M_C$
carries additional regularity relative to $\M$. The numerical experiments use a Mat\'ern operator for
this purpose. It is chosen after training and is not used to construct the
derivative-supervision directions, although its hyperparameters lie within the
range of the Mat\'ern family used to generate training parameters.

Working in this stronger space entails the regularity condition
$m^\dagger-m_0\in\M_C$ of Assumption~\ref{ass:preconditioner}. For a truth
outside $m_0+\M_C$, the iteration \eqref{eq:preconditioned-surrogate-lw} is
still well defined, but Theorem~\ref{thm:convergence} does not
apply. Remark~\ref{rem:preconditioner} describes the condition for
operators of Mat\'ern type.

\subsection{Which derivative errors matter under preconditioning}
\label{sec:weighted-derivative-error}

Changing the parameter-space geometry changes the derivative approximation
that enters the Landweber iteration.  With $B_\theta(m)$ as in \eqref{eq:derivative-error-operator}, we write
$\varepsilon_{D,C}(m):=\|B_\theta(m)C^{1/2}\|_{\Lop(\M,\Y)}$ for the
pointwise weighted error, whose supremum over $B_\rho^C(m^\dagger)$ is
$\varepsilon_{D,C}$ in \eqref{eq:accuracy}; it satisfies
$\varepsilon_{D,C}(m)=\|C^{1/2}B_\theta(m)^*\|_{\Lop(\Y,\M)}$.
For any residual $R\in\Y$,
\begin{equation}
    \|CB_\theta(m)^*R\|_{\M_C}
    \leq
    \varepsilon_{D,C}(m)\|R\|_{\Y}.
    \label{eq:weighted-adjoint-error}
\end{equation}
The derivative error that matters for \eqref{eq:preconditioned-surrogate-lw}
is therefore the weighted error $\varepsilon_{D,C}(m)$, whatever
the training. Smoother direction spaces also appear in the DIFNO theory of
\citet{Yao2025}, but there they are part of the approximation result, whereas
here the stronger space is chosen when the trained surrogate is used for inversion.

The weighted error connects the inversion geometry with finite directional
supervision.  Let $P_r$ be the orthogonal projection onto the supervised
subspace $H_r=\operatorname{span}\{h_1,\ldots,h_r\}$.  Decomposing with respect
to this subspace gives
\begin{equation}
    \|B_\theta(m)C^{1/2}\|
    \leq
    \|B_\theta(m)P_r\|\,\|P_rC^{1/2}\|
    +\|B_\theta(m)(I-P_r)\|\,\|(I-P_r)C^{1/2}\|.
    \label{eq:supervision-preconditioner-bound}
\end{equation}
The first term involves only the supervised directions and is therefore
related to the error measured by the directional training objective.  The second term
contains derivative components not controlled by the prescribed directions,
weighted by the extent to which the chosen parameter-space geometry retains
them.

For a smoothing $C$, rough parameter directions are attenuated by $C^{1/2}$,
which makes the iteration less sensitive to poorly approximated derivative
components. It does not recover information that was never learned. The same
operator acts on the exact gradient and attenuates genuine fine-scale features
as well as surrogate error, so smoothing can improve robustness to continued
iteration while slowing the recovery of fine-scale structure.

Forward error remains a separate contribution. With $b=\F_\theta(m)-\F(m)$,
the learned and exact preconditioned search directions differ by
$C(D\F_\theta(m)-D\F(m))^*R_\theta+C\,D\F(m)^*b$, where
$R_\theta=\F_\theta(m)-y^\delta$; in the norm of $\M_C$ the first term is at
most $\varepsilon_{D,C}(m)\|R_\theta\|_\Y$, whereas the second is at most
$\|D\F(m)C^{1/2}\|\,\|b\|_\Y$ and is controlled only by the forward error. Changing the geometry can therefore reduce the influence of
derivative error, but it does not remove the forward error and does not
guarantee a more accurate search direction in relative terms.

\subsection{Relation to the regularisation analysis}
\label{sec:preconditioned-regularisation}

For compact $C$, Lemma~\ref{lem:descent} and Theorem~\ref{thm:convergence}
show that the effective error level in the discrepancy principle is $\Delta$
in \eqref{eq:effective_noiselevel}, with $\gamma=\rho\,\varepsilon_{D,C}$
under \eqref{eq:tcc} and $\gamma=0$ under \eqref{eq:tcc-surrogate}. The
derivative error therefore enters only through its weighted form, and not at
all when the surrogate itself satisfies the cone condition. The scale of $C$
plays no role: replacing $C$ by $sC$ multiplies $\varepsilon_{D,C}$ by
$s^{1/2}$ and the radius $\rho$ in $\M_C$ by $s^{-1/2}$, so that only the
relative weighting of parameter directions matters. A smoothing $C$ reduces
$\varepsilon_{D,C}$ relative to the unweighted derivative error, but it also
increases $\rho$, so the product $\rho\,\varepsilon_{D,C}$ can
decrease or increase.

The existence of surrogates with small $\varepsilon_{\mathrm{op}}$ and
$\varepsilon_{D,C}$ on $B_\rho^C(m^\dagger)$, noted after
\eqref{eq:accuracy}, is an existence result. Minimising the derivative term
in \eqref{eq:empirical-reduced-dino-loss} on finitely many parameters and
directions does not establish a uniform bound on $\varepsilon_{D,C}$. The numerical experiments keep the
derivative-supervision protocol fixed while varying the parameter-space
geometry used in the inversion. This separates the information used to construct the learned operator from
the geometry in which it is used.

\section{Numerical experiments}
\label{sec:numerical-interpretation}

\paragraph{Implementation.}
All finite-element computations, in this section and in
Section~\ref{sec:darcy}, use FEniCSx/DOLFINx~0.10.0 with PETSc~3.25.0, and
exact derivatives are computed from the discrete tangent and adjoint
equations of the same discretisation. The neural operators are implemented
in PyTorch~2.4.0 with CUDA~12.4. GPU computations use the default,
non-deterministic reductions of PyTorch, so repeated surrogate runs are not
reproducible bit for bit, and rounding differences can grow along an
iteration. Finite-element solves run on CPU cores, eight per exact
inversion. For the reaction benchmark, networks are trained on one NVIDIA
A100 GPU and surrogate inversions run on one NVIDIA L40 GPU. All runs were
performed on the Sulis Tier~2 HPC platform.

\subsection{Reaction benchmark: forward model and discretisation}
\label{sec:reaction-benchmark}

The principal benchmark is the recovery of a spatially varying reaction
coefficient in an elliptic PDE, simple enough that forward, derivative and
adjoint computations can be evaluated accurately by finite elements, while
retaining the nonlinear parameter-to-state structure relevant to the
analysis. The unknown is a reaction coefficient $m$ on $\Omega=(0,1)^2$, with
state $p$ satisfying
\begin{equation}
    -\Delta p + m p = 1
    \quad \text{in }\Omega,
    \qquad
    p=0
    \quad \text{on }\partial\Omega ,
    \label{eq:reaction-forward}
\end{equation}
and $\F(m)=p(m)$ denotes the parameter-to-state map, with parameter and data
spaces $L^2(\Omega)$. The admissible coefficients satisfy $0\leq m\leq 200$,
so that the weak problem is coercive. For a perturbation $h$, the
directional derivative $v=D\F(m)h$ satisfies
\begin{equation}
    -\Delta v + mv = -hp
    \quad \text{in }\Omega,
    \qquad
    v=0
    \quad \text{on }\partial\Omega ,
    \label{eq:reaction-derivative}
\end{equation}
and, for $z\in L^2(\Omega)$, the adjoint action is $D\F(m)^*z=-pw$, where
\begin{equation}
    -\Delta w + mw = z
    \quad \text{in }\Omega,
    \qquad
    w=0
    \quad \text{on }\partial\Omega ,
    \label{eq:reaction-adjoint}
\end{equation}
so that both the derivative labels used in training and the exact Landweber
directions require only elliptic solves.

\paragraph{Discretisation.}
The surrogate training and inverse calculations use continuous
piecewise-linear finite elements on a uniform triangulation obtained from
$119\times119$ square cells, giving a $120\times120$ nodal representation of
the coefficient. The state is observed at the interior mesh vertices. To
keep the notation light, we write $\F$ and $D\F$ also for the resulting
discrete forward map and its derivative; all computations below use the
discrete model.

\subsection{Training measure}
\label{sec:reaction-training-measure}

The measure $\mu$ of Section~\ref{sec:directional-training} is a
heterogeneous family of stationary Mat\'ern fields: training coefficients
are $m=a+\sigma z_{\nu,\ell}$, where $z_{\nu,\ell}$ has zero mean and unit
marginal variance, $a\sim\mathcal U(95,105)$, $\sigma\sim\mathcal U(10,30)$,
$\nu\in\{1,1.5,2,2.5\}$ and $\ell\in\{0.075,0.10,0.15,0.20,0.25\}$. When
necessary, the fluctuation about $a$ is rescaled to enforce
$0\leq m\leq200$, so $\mu$ is a non-Gaussian mixture over amplitudes,
regularities and correlation lengths. It is the measure in the risk
\eqref{eq:population-dino-risk}; the inversion itself uses no prior. We draw
$10\,000$ coefficients for training, $2\,500$ for validation and $2\,500$ for
final testing, independently, and a further independent ensemble of
$10\,000$ fields to construct the derivative directions.

\subsection{Learning the forward map and derivative actions}
\label{sec:reaction-learning}

We compare four learned forward operators, each mapping the nodal
coefficient to the nodal state on the $120\times120$ grid. All are trained
with AdamW for $2500$ epochs (learning rate $10^{-3}$ for the FNO and
$10^{-4}$ otherwise), with the learning rate reduced on plateaus of the
validation forward error. Their parameter counts are listed in
Table~\ref{tab:headline_emulator}.
\begin{itemize}
 \item \textbf{FNO.} A standard FNO \citep{Li2021} whose lifting receives the
 normalised coefficient and the coordinates $(x,y)$, followed by four
 Fourier layers of width $24$ with $10\times10$ retained modes, padding $8$,
 GELU activations and a projection of width $128$.
 \item \textbf{DeepONet.} A standard DeepONet \citep{Lu2021}. The branch is
 a convolutional network with GroupNorm, ReLU and max pooling acting on the
 normalised coefficient and the coordinates, and returns $400$ coefficients.
 The trunk has six GELU layers of width $400$ acting on Fourier features of
 $(x,y)$ with frequencies $1,2,\ldots,32$. The output is the scaled inner
 product of branch and trunk plus a bias.
 \item \textbf{GenDeepONet.} A generalised DeepONet \eqref{eq:gdon}. The
 encoder subtracts the mean field and returns the standardised coordinates
 of $m$ in the leading $400$ input PCA modes, and the decoder combines the
 leading $300$ output POD modes and adds the mean state. Both bases use
 trapezoidal quadrature weights and are computed from the $10\,000$ training
 fields before training. The network $f_\theta$ is a residual multilayer
 perceptron of width $1024$ with four GELU residual blocks.
 \item \textbf{GenFNO.} A generalised Fourier operator \eqref{eq:gfno}. The
 encoder applies a discrete Fourier transform to the normalised coefficient
 and keeps a rectangular block of low frequencies, $r=1024$ real
 coordinates; $f_\theta$ is a residual multilayer perceptron of width $512$
 with four GELU residual blocks, returning $s=2304$ real Fourier
 coefficients from which the decoder synthesises the state. The Dirichlet
 condition is not imposed but learned from the data, and the inversion uses
 only interior observations.
\end{itemize}
All derivatives of the surrogates are exact derivatives of the complete
networks: Jacobian--vector products in the derivative loss and
vector--Jacobian products, with the residual held fixed, in the Landweber
update. They are computed by automatic differentiation of the nonlinear part,
with the fixed linear encoders and decoders either included in the
differentiation or applied analytically; no separate derivative network is
trained. GenDeepONet and GenFNO are instances of the encoder--network--decoder
class $\Phi f_\theta\Psi$ of \eqref{eq:rb}, with fixed $\Psi$ and $\Phi$
given by PCA/POD and by Fourier coordinates, respectively (the mean fields and
coordinate scalings are affine and can be absorbed as in
Section~\ref{sec:generalised-deeponet}). The results of
Theorems~\ref{thm:compact-directions} and~\ref{thm:operator-approx}
therefore apply to their classes through
Corollary~\ref{cor:architecture-transfer}, since GELU is $C^1$ and
nonpolynomial, although they concern the classes, not the sizes used here,
and for GenFNO the Fourier coordinates on the square lie outside the
periodic setting of the corollary. The standard FNO and DeepONet are not of
this form: the FNO couples frequencies only through pointwise nonlinearities
(Remark~\ref{rem:fno-frequency}), and the DeepONet learns its output
functions through the trunk (Remark~\ref{rem:trained-trunk}) and, with its
ReLU branch, is only piecewise smooth in $m$. The purpose of the comparison
is not to rank architectures, but to examine whether learned operators with
different representations can supply sufficiently accurate values and
derivatives for iterative regularisation.

\paragraph{Derivative directions.}
The $r=16$ training directions of Section~\ref{sec:directional-training} are
built from the independent ensemble. We separate the
spatial mean $\langle m\rangle_\Omega$ of each field from its fluctuation
$m'=m-\langle m\rangle_\Omega\mathbf 1$, apply PCA to the fluctuations, and
take
\begin{equation}
    h_1=\sigma_{\mathrm{mean}}\,\mathbf 1,
    \qquad
    h_{j+1}=\sqrt{\kappa_j}\,\psi_j,\quad j=1,\ldots,r-1,
    \label{eq:reaction-derivative-directions}
\end{equation}
where $\sigma_{\mathrm{mean}}$ is the standard deviation of the spatial mean
and $(\kappa_j,\psi_j)$ are the leading empirical variances and modes of the
fluctuations. Because the mean is split off first, these directions
differ from the $r$ leading eigenvectors of the empirical covariance of $m$.
The derivative loss therefore weights derivative errors by
$\widehat\Gamma_r=\sum_{j=1}^{r}h_j\otimes h_j$, in the notation of
Section~\ref{sec:directional-training}, and not by a truncated covariance
$\Gamma_{\mu,r}$.

\paragraph{Training and evaluation.}
For every training coefficient $m_i$, the finite-element model supplies
$\F(m_i)$ and the derivative actions $D\F(m_i)h_j$, $j=1,\ldots,r$, from
\eqref{eq:reaction-derivative}, and the same data are given to all four
architectures. During derivative-informed training, one direction is drawn
uniformly at random for each occurrence of a training coefficient, as in
Section~\ref{sec:directional-training}. Since each architecture computes its
loss in its own normalised or reduced coordinates, all comparisons use
common physical-space metrics. On the independent test set, with $N_F=2500$
fields, a fixed subset $I_J$ of $256$ fields for derivatives and the
full-grid trapezoidal weights $W_p$, we report
\begin{equation}
    E_F^2=\frac{\sum_{i=1}^{N_F}\|\F_\theta(m_i)-\F(m_i)\|_{W_p}^2}
               {\sum_{i=1}^{N_F}\|\F(m_i)\|_{W_p}^2},
    \qquad
    E_J^2=\frac{\sum_{i\in I_J}\sum_{j=1}^{r}\|B_\theta(m_i)h_j\|_{W_p}^2}
               {\sum_{i\in I_J}\sum_{j=1}^{r}\|D\F(m_i)h_j\|_{W_p}^2},
    \label{eq:reaction-EF-EJ}
\end{equation}
where $B_\theta$ is the derivative error \eqref{eq:derivative-error-operator}.
$E_J$ measures generalisation to new coefficient fields along the training
directions, that is, in the $\widehat\Gamma_r$ geometry; it is not the full
Hilbert--Schmidt or operator-norm error, nor a bound along an inverse
trajectory.

For each architecture we considered
$\lambda_J\in\{0,\,0.1,\,1,\,10,\,100,\,300,\,1000,\,3000\}$ and selected the
model and checkpoint on the validation set only, by the score
$S_{\mathrm{val}}=[(E_F^{\mathrm{val}})^2+(E_J^{\mathrm{val}})^2]^{1/2}$.
Table~\ref{tab:headline_emulator} gives the selected weights and epochs and
the test errors: all four surrogates have forward errors below $0.5\%$ and
derivative errors below $0.8\%$. For each architecture we also retain the
model trained without the derivative loss ($\lambda_J=0$), whose checkpoint
is selected by the same score; this gives it derivative information its
training did not use, so the comparison in
Section~\ref{sec:reaction-derivative-training}, if anything, favours it.

\begin{table}[t]
    \centering
    \caption{Independent-test accuracy of the four validation-selected
    surrogates.}
    \label{tab:headline_emulator}
    {\small\begin{tabular}{lrrrrrr}
\toprule
Architecture & $\lambda_J$ & Parameters & Epoch & $E_F$ (\%) & $E_J$ (\%) & $S$ (\%) \\
\midrule
FNO & 300 & 466,625 & 2300 & 0.091 & 0.157 & 0.181 \\
DeepONet & 100 & 15,535,697 & 2500 & 0.166 & 0.769 & 0.786 \\
GenFNO & 100 & 3,808,000 & 2500 & 0.497 & 0.478 & 0.690 \\
GenDeepONet & 1000 & 9,114,924 & 2500 & 0.213 & 0.424 & 0.475 \\
\bottomrule
\end{tabular}
}
\end{table}

\subsection{Inversion setup}
\label{sec:reaction-inversion-setup}

\paragraph{Inverse truths and data.}
We consider two prescribed inverse truths. T1 is an independent draw from the
training family, with $\nu=3/2$ and $\ell=0.075$, and provides an
in-distribution problem. T2 is a deterministic coefficient with background
value $100$, a rectangular inclusion of value $150$ and two circular
inclusions of values $20$ and $180$; it lies in the admissible range but has
sharp geometric structure absent from the training family. The observations
are computed on a finer $952\times952$-cell mesh and restricted to the
observation grid, so that the data are not generated with the
discretisation used for inversion. The comparisons below are noiseless until
Section~\ref{sec:reaction-noise}.

\paragraph{Landweber implementation.}
The baseline inverse calculations use the Mat\'ern parameter-space geometry
introduced in Section~\ref{sec:preconditioned-landweber}, with
$\nu_C=3/2$, $\ell_C=0.1$ and $\sigma_C^2=1$, and start from $m_0\equiv100$.
Let $W_y$ denote the quadrature weights in the discrete observation norm,
$R_k=\F_\theta(m_k)-y_h^\delta$ the residual and
$J_k=D\F_\theta(m_k)$ the derivative of the learned map (or $\F$ and $D\F$
for the exact iteration). We form
$q_k=J_k^\top W_yR_k$ and $d_k=C_hq_k$, where $C_h$ is the discrete action
of the chosen parameter-space operator, and use the adaptive step
\begin{equation}
    \alpha_k
    =
    \min\left\{
        0.95\,
        \frac{q_k^\top d_k}
             {\|J_kd_k\|_{W_y}^2},
        \frac{5}{\|d_k\|_\infty}
    \right\}.
    \label{eq:reaction-adaptive-step}
\end{equation}
The proposed update $m_k-\alpha_kd_k$ is projected componentwise onto
$[0,200]$, without backtracking, and the exact-FEM reference uses the same
protocol with finite-element evaluations. The first term in
\eqref{eq:reaction-adaptive-step} is $0.95$ times the steepest-descent step in
the geometry induced by $C_h$, which the descent argument of
Lemma~\ref{lem:descent} covers (Remark~\ref{rem:step-sizes} on variable step
sizes). The $\ell^\infty$ cap and the projection are not covered by the
analysis; see \citet{PianaBertero1997} for projected Landweber iteration in
the linear case.

The matrix $C_h$ is the nodal Mat\'ern covariance matrix
$(C_h)_{ij}=k(\|x_i-x_j\|)$ for grid nodes $x_i,x_j$, with the Mat\'ern
kernel \citep[Section~4.2]{RasmussenWilliams2006}
\begin{equation}
    k(r)=\sigma_C^2\,\frac{2^{1-\nu_C}}{\Gamma(\nu_C)}
    \Bigl(\frac{\sqrt{2\nu_C}\,r}{\ell_C}\Bigr)^{\nu_C}
    K_{\nu_C}\Bigl(\frac{\sqrt{2\nu_C}\,r}{\ell_C}\Bigr),
    \label{eq:matern-kernel}
\end{equation}
where $K_{\nu}$ is the modified Bessel function of the second kind. It is
constructed from the grid geometry alone, independently of the training
data, and applied exactly, up to rounding, by zero-padded FFT convolution,
so that opposite boundaries are not identified. It discretises the
covariance operator $(Cf)(x)=\int_\Omega k(\|x-y\|)f(y)\,\mathrm dy$ on
$L^2(\Omega)$, the quadrature factor being absorbed by the adaptive step, so
by Remark~\ref{rem:preconditioner} the baseline geometry has
$\M_C=H^{5/2}(\Omega)$. Neither truth satisfies $m^\dagger-m_0\in\M_C$: T2
is discontinuous, and T1, a Mat\'ern sample with the same smoothness
$\nu=3/2$ as $C$, lies outside its Cameron--Martin space almost surely. The
baseline runs therefore lie outside the hypotheses of
Theorem~\ref{thm:convergence}; Section~\ref{sec:reaction-low-smoothness}
repeats the noiseless experiments with geometries ($\nu_C<1/2$) for which
the regularity condition holds for T1.

With $m_h^\dagger$ the target coefficient on the inversion grid and $W_m$
the corresponding quadrature weights, we report the relative coefficient
error and the oracle iteration
\begin{equation}
    e_m(k)=\frac{\|m_k-m_h^\dagger\|_{W_m}}{\|m_h^\dagger\|_{W_m}},
    \qquad
    k_{\mathrm{or}}=\argmin_{0\leq k\leq10\,000}e_m(k).
    \label{eq:reaction-coefficient-error}
\end{equation}
The oracle uses the known truth; it serves only to compare the forward
models and is not a stopping rule. Many oracle iterations coincide with the
end of the $10\,000$-iteration budget, and the corresponding errors are then
finite-budget values rather than minima.

\subsection{Noiseless reconstruction with learned operators}
\label{sec:reaction-noiseless}

\begin{table}[tbp]
    \centering
    \caption{Noiseless reaction-coefficient inversion for the held-out
    in-distribution truth T1 and the structured truth T2.  Here $e_m^\star=e_m(k_{\mathrm{or}})$ and
    $e_m^K$ is the error at the final iterate $K=10\,000$. The residuals are
    relative data misfits at the oracle iterate,
    $100\,\|G(m_{k_{\mathrm{or}}})-y_h\|_{W_y}/\|y_h\|_{W_y}$, with
    $G=\F_\theta$ for the surrogate residual and $G=\F$ for the FEM
    residual. The latter re-evaluates the surrogate reconstruction with the
    finite-element model. For exact FEM only the FEM residual applies.}
    \label{tab:reaction_noiseless_headline}
    \resizebox{\linewidth}{!}{%
\begin{tabular}{llrrrrr}
\toprule
Truth & Method & $k_{\mathrm{or}}$ & $e_m^\star$ (\%) & $e_m^{K}$ (\%) & Surrogate residual (\%) & FEM residual (\%) \\
\midrule
T1 & Exact FEM & 10000 & 4.118 & 4.118 & -- & 0.035 \\
T1 & FNO & 10000 & 6.639 & 6.639 & 0.149 & 0.159 \\
T1 & DeepONet & 10000 & 7.651 & 7.651 & 0.235 & 0.236 \\
T1 & GenFNO & 10000 & 5.825 & 5.825 & 0.462 & 0.103 \\
T1 & GenDeepONet & 9143 & 7.105 & 7.106 & 0.178 & 0.228 \\
\midrule
T2 & Exact FEM & 10000 & 6.936 & 6.936 & -- & 0.068 \\
T2 & FNO & 10000 & 8.858 & 8.858 & 0.213 & 0.214 \\
T2 & DeepONet & 10000 & 9.703 & 9.703 & 0.334 & 0.327 \\
T2 & GenFNO & 8859 & 8.574 & 8.577 & 0.462 & 0.371 \\
T2 & GenDeepONet & 5523 & 9.405 & 9.454 & 0.219 & 0.579 \\
\bottomrule
\end{tabular}
}

\end{table}

Table~\ref{tab:reaction_noiseless_headline} summarises the noiseless
reconstructions for T1 and T2 with the four selected surrogates and the
baseline geometry.  The exact-FEM iteration reaches oracle errors of $4.12\%$ on T1 and $6.94\%$ on T2.  The corresponding
surrogate oracle errors range from $5.83\%$ to $7.65\%$ on T1 and from
$8.57\%$ to $9.70\%$ on T2.  All four learned operators therefore support
meaningful reconstruction on both the in-distribution and structured truths,
although with an accuracy loss relative to the finite-element reference.

The last two columns compare the misfit seen by each surrogate with the
misfit of the same reconstruction under the finite-element model. For FNO and
DeepONet the two agree closely. For GenFNO on T1 the surrogate reports
$0.46\%$ against a finite-element misfit of $0.10\%$, whereas for
GenDeepONet on T2 it reports $0.22\%$ against $0.58\%$. The residual of the
surrogate therefore need not reflect how well its reconstruction fits the
data: it can overstate the true misfit, as for GenFNO on T1, or understate
it, as for GenDeepONet on T2.

\begin{figure}[!t]
    \centering
    \includegraphics{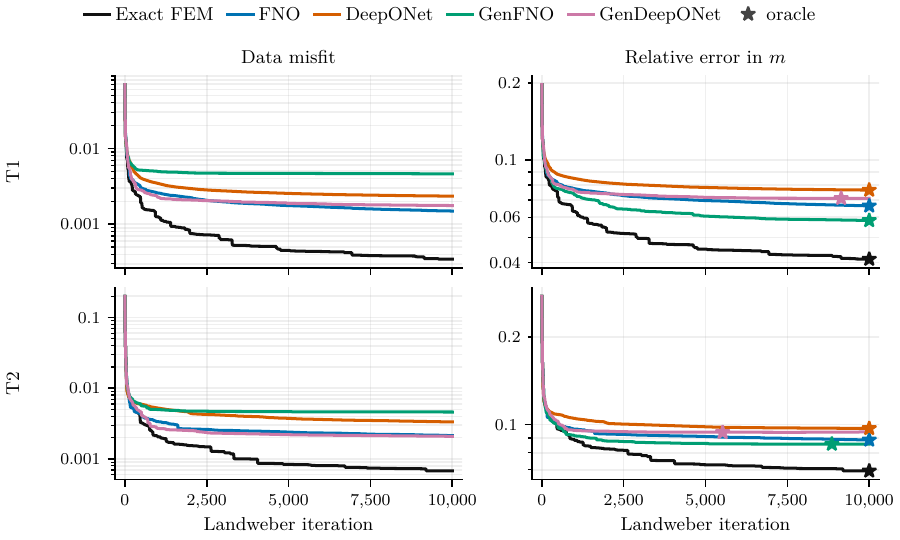}

    \caption{Noiseless data misfit and relative coefficient error for T1
    (top) and T2 (bottom), for the exact finite-element iteration and the
    four selected learned operators.
    Stars mark the smallest coefficient error over the $10\,000$ iterations. Each data misfit is evaluated with the forward
    model used in that iteration.}
    \label{fig:reaction-noiseless-convergence}
\end{figure}

Figure~\ref{fig:reaction-noiseless-convergence} shows the complete
noiseless trajectories. With exact data the exact iteration keeps reducing
the coefficient error over all $10\,000$ steps, as expected for noise-free
data, and so do most surrogate iterations. The exceptions are GenFNO and
GenDeepONet on T2, whose error is smallest before the end of the run
(Table~\ref{tab:reaction_noiseless_headline}) and then grows slightly while
their residual keeps decreasing. Here the surrogate error plays the role of
data noise, and this is the onset of the semiconvergence studied in
Section~\ref{sec:reaction-geometry}.

\begin{figure}[tbp]
    \centering
    \includegraphics[width=\linewidth]
    {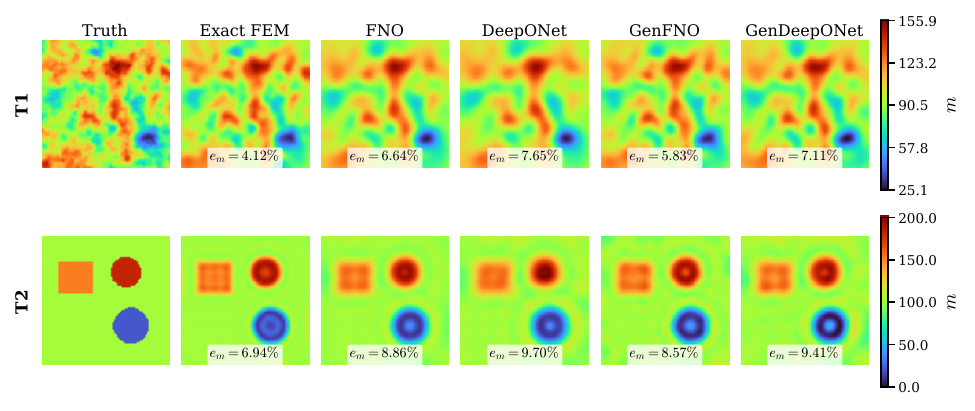}
    \caption{Noiseless T1 (top) and T2 (bottom) truths and oracle
    reconstructions $m_{k_{\mathrm{or}}}$ with the baseline Mat\'ern geometry.
    Each row has its own colour scale, from the minimum to the maximum over
    its six fields. The inclusions of T2 are not drawn from the random-field
    morphology used for training.}
    \label{fig:oracle_fields}
\end{figure}

Figure~\ref{fig:oracle_fields} shows the oracle reconstructions: every
method reproduces the large-scale variation of T1 and, although the training
distribution contains no such shapes, the locations and relative contrasts of
the three inclusions of T2. Even the exact-FEM reconstruction smooths the
interfaces of T2, so the loss of sharpness is due partly to the baseline
geometry, which Section~\ref{sec:reaction-geometry} varies.

\subsection{Effect of derivative-informed training}
\label{sec:reaction-derivative-training}

To isolate the effect of derivative supervision, we compare each selected
surrogate with the architecture-matched model trained without the derivative
loss ($\lambda_J=0$), with everything else unchanged. Derivative supervision
reduces the test derivative error $E_J$ by $45\%$ for FNO, $72\%$ for
DeepONet, $63\%$ for GenFNO and $94\%$ for GenDeepONet, but hardly changes
the best reconstruction: the oracle errors change by between $-1.1$ and
$+0.4$ percentage points, with mixed sign, which is of the size to be expected
between training seeds for models trained once. Its clearer effect is on the
late iterations. Without the derivative loss, GenFNO and GenDeepONet reach
their smallest errors early and then deteriorate, on T2 at iterations $7185$
and $1286$; with it these move to $8859$ and $5523$, and on T1 from $7587$ and
$5601$ to $10\,000$ and $9143$. FNO and DeepONet reach their smallest errors at
the end of the run in either case. Derivative accuracy thus delays
semiconvergence, as the level $\Deff$ predicts, but does not by itself
determine the best reconstruction; more generally, smaller $E_F$ or $E_J$ in
Table~\ref{tab:headline_emulator} does not imply a smaller coefficient error in
Table~\ref{tab:reaction_noiseless_headline}. We examine this distinction in
the next subsection.

\subsection{From directional derivative accuracy to the Landweber direction}
\label{sec:reaction-derivative-geometry}

The metric $E_J$ measures derivative accuracy along the training directions,
in the geometry of $\widehat\Gamma_r$, whereas the baseline iteration depends
on the weighted adjoint action $C\,D\F_\theta(m_k)^*R_k$; as discussed in
Section~\ref{sec:training-preconditioning}, accuracy in the first need not
imply accuracy in the second.

To examine this distinction without confounding it with
architecture-dependent inverse trajectories, we evaluate the derivative
diagnostics at two common coefficient fields, the initial point
$m_0\equiv100$ and the T2 truth.  At such a field $m$, let $R=\F_\theta(m)-y_h$ be the residual of the
surrogate with respect to the noiseless T2 data $y_h$, and let
$q_\theta=D\F_\theta(m)^\top W_yR$ and $q=D\F(m)^\top W_yR$ be the learned and exact adjoint actions, that is, the vector--Jacobian
products (VJPs) for this residual. Using the same $R$ in both isolates the
error of the derivative from that of the forward map. With $B_\theta$
as in \eqref{eq:derivative-error-operator}, we compare three relative errors,
\[
    \Biggl[\frac{\sum_{j}\|B_\theta(m)h_j\|_{W_y}^2}
                 {\sum_{j}\|D\F(m)h_j\|_{W_y}^2}\Biggr]^{1/2},
    \qquad
    \frac{\|q_\theta-q\|_2}{\|q\|_2},
    \qquad
    \frac{\bigl[(q_\theta-q)^\top C_h(q_\theta-q)\bigr]^{1/2}}
         {\bigl[q^\top C_hq\bigr]^{1/2}}.
\]
The first is the pooled error of the Jacobian--vector products (JVPs) along
the $r$ training directions $h_j$ of \eqref{eq:reaction-derivative-directions}.
The second is the error of the VJP in the Euclidean norm, which governs an
unpreconditioned Landweber step. The third is the error of the
preconditioned directions $C_hq_\theta$ and $C_hq$ in the norm of $\M_C$, in
which Lemma~\ref{lem:descent} measures the iteration. With
the discrete inner product $\langle u,v\rangle_{\M_C}=u^\top C_h^{-1}v$, one
has $\|C_hv\|_{\M_C}^2=v^\top C_hv$. The residual $R$ depends on the
surrogate, so the VJP errors of different architectures use slightly
different residuals.

\begin{figure}[tbp]
    \centering
        \includegraphics{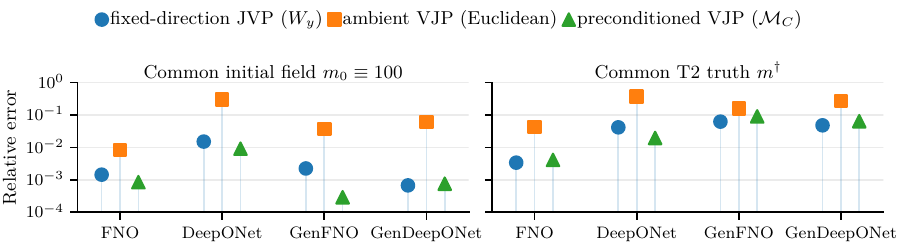}
    
    \caption{The three relative derivative errors defined in the text, at the
    initial field $m_0\equiv100$ and at the T2 truth: pooled JVP error along
    the training directions, Euclidean (ambient) VJP error, and
    preconditioned VJP error in the norm of $\M_C$.}
    \label{fig:derivative_geometry}
\end{figure}

Figure~\ref{fig:derivative_geometry} shows that these notions of derivative
accuracy can differ substantially.  At the initial field, for example,
DeepONet has a pooled fixed-direction JVP error of $1.52\%$, but its ambient
same-residual VJP error is $30.08\%$.  After application of the Mat\'ern
operator, the corresponding VJP error, measured in $\M_C$, is $1.30\%$.  The relative separation is even larger for GenDeepONet, whose fixed-direction error at $m_0$ is only $0.067\%$ against an ambient VJP error of $6.02\%$.

The same qualitative distinction remains at the T2 truth, although its
magnitude is model dependent.  In every case the ambient VJP error is larger
than the corresponding fixed-direction error, and applying the Mat\'ern
geometry reduces the ambient adjoint error.  The preconditioned error at the T2 truth nevertheless remains appreciable for GenFNO and
GenDeepONet.

This helps explain why the reductions in $E_J$ reported in
Section~\ref{sec:reaction-derivative-training} do not translate into better
reconstructions: the training loss probes a few fixed input directions,
whereas the iteration uses the adjoint action for its current residual,
weighted by $C$.

\subsection{Preconditioning and semiconvergence}
\label{sec:reaction-geometry}

Section~\ref{sec:landweber} motivates the preconditioner by the fact that the
learned derivative is unreliable in rough directions, which are essentially
absent from the training data. To test this, we run the same noiseless T2
inversions with the identity in place of the Mat\'ern geometry, keeping the
network weights, data and initial field. This setting lies outside the theory
of Section~\ref{sec:landweber}. The Mat\'ern results are those of
Table~\ref{tab:reaction_noiseless_headline}, and the smallest errors with the
identity are also listed in Table~\ref{tab:preconditioner_low_smoothness}.

\begin{figure}[tbp]
    \centering
    \includegraphics[width=\linewidth]
    {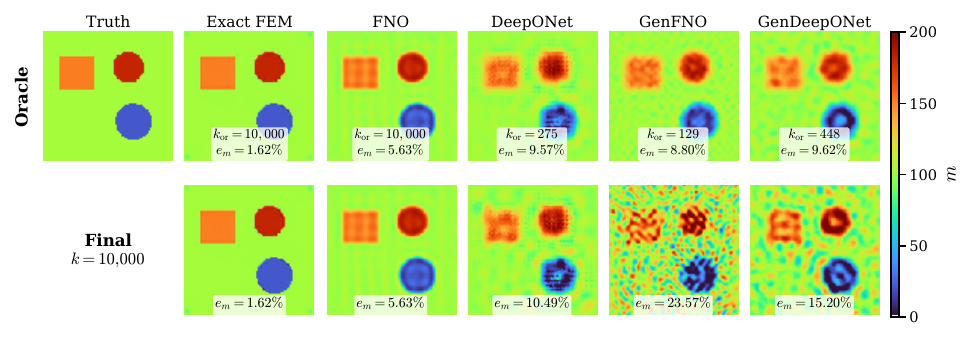}
    \caption{Noiseless T2 reconstructions without preconditioning. Top: each
    method at its oracle iterate. Bottom: the final iterate, $k=10\,000$.
    Continued iteration produces oscillatory artefacts and projection onto the
    bounds for DeepONet, GenFNO and GenDeepONet.}
    \label{fig:identity_fields}
\end{figure}

Without preconditioning, three of the four surrogates semiconverge although
the data are exact, whereas the exact iteration and the FNO reach their
smallest errors at the end of the run. DeepONet reaches its minimum of $9.57\%$ after $275$
iterations and ends at $10.49\%$. GenFNO reaches $8.80\%$ at iteration $129$
and ends at $23.57\%$, and GenDeepONet reaches $9.62\%$ at iteration $448$
and ends at $15.20\%$. With the Mat\'ern geometry the same models stay within
$1\%$ of their minima over the budget. The loss of accuracy comes with
spurious structure (Figure~\ref{fig:identity_fields}): fine oscillatory
texture already at the oracle iterates, and by the end of the run
high-frequency speckle, ring-like errors around the inclusions and values
pushed to the bounds $0$ and $200$. The exact-FEM control shows that the
problem itself does not require the smoothing. With the identity, exact
Landweber reaches $1.62\%$ on T2, against $6.94\%$ with the Mat\'ern geometry,
and the FNO likewise improves from $8.86\%$ to $5.63\%$. We also decreased $\ell_C$ from $0.1$ to $0.005$, which moves the geometry
gradually from the baseline towards the identity (results not shown). Along
this path the best FNO error decreases steadily, reaching $0.64$ times its
baseline value under the identity, while those of the other surrogates change little. This is consistent
with Figure~\ref{fig:derivative_geometry}. The FNO has much the smallest
unpreconditioned VJP error, so its unfiltered direction remains accurate,
whereas for the other surrogates the components removed by the smoothing are
dominated by derivative error.

The instability without preconditioning therefore comes from the surrogates,
not from the inverse problem. As anticipated in
Section~\ref{sec:weighted-derivative-error}, the Mat\'ern geometry does not
make the learned derivative more accurate but damps its poorly learned
fine-scale components, together with genuine fine-scale information, which is
why it improves stability while slowing the recovery of the sharp features of
T2.

\paragraph{Surrogate error as a source of semiconvergence.}
\label{sec:reaction-capacity}
In Lemma~\ref{lem:descent} the surrogate errors enter the level
$\Delta=\delta+\epsop+\rho\,\varepsilon_{D,C}/(1+\eta)$ in the same way as the
noise.
Even with exact data, descent is guaranteed only while the residual exceeds
$\tau\Deff$, and beyond that point an inaccurate surrogate can drive the
iteration away from the truth, as noise would. To test this we use a
capacity study of the FNO, in which the width and number of Fourier modes are
varied while the training data, the derivative weight $\lambda_J=10$ and the
selection rule are fixed (Table~\ref{tab:fno_capacity}), and a smaller study
in which the number of training fields is also reduced
(Table~\ref{tab:fno_training_size}).

\begin{table}[!t]
    \centering
    \caption{FNO capacity study at $\lambda_J=10$, noiseless T2 data. For
    each capacity the same checkpoint is used with both geometries. A
    final-to-minimum ratio above one indicates semiconvergence within the
    $10\,000$-iteration budget.}
    \label{tab:fno_capacity}
        {\small\begin{tabular}{lrr@{\hspace{10pt}}rrr@{\hspace{10pt}}rrr}
\toprule
 & & & \multicolumn{3}{c}{Mat\'ern, $\nu_C=3/2$, $\ell_C=0.1$} & \multicolumn{3}{c}{Identity} \\
\cmidrule(lr){4-6}\cmidrule(lr){7-9}
Width/modes & $E_F$ (\%) & $E_J$ (\%) & $e_m^\star$ (\%) & $k_{\rm or}$ & Final/min & $e_m^\star$ (\%) & $k_{\rm or}$ & Final/min \\
\midrule
8/4 & 1.357 & 2.487 & 11.79 & 248 & 1.21 & 9.68 & 87 & 2.13 \\
12/6 & 0.905 & 0.956 & 10.65 & 551 & 1.23 & 8.92 & 116 & 1.79 \\
16/8 & 0.264 & 0.383 & 9.60 & 9989 & 1.00 & 7.22 & 1115 & 1.05 \\
24/10 & 0.115 & 0.224 & 8.88 & 10000 & 1.00 & 5.52 & 10000 & 1.00 \\
32/12 & 0.018 & 0.118 & 8.84 & 10000 & 1.00 & 5.13 & 10000 & 1.00 \\
40/16 & 0.014 & 0.113 & 8.78 & 10000 & 1.00 & 5.41 & 10000 & 1.00 \\
\bottomrule
\end{tabular}
}

    \vspace{1.2em}
    \caption{FNO models trained on fewer fields, at $\lambda_J=10$, with the
    matching $10\,000$-field models for reference. Noiseless inversion with
    the baseline Mat\'ern geometry. Errors $E_F$ and $E_J$ are validation
    errors. A dash means that the inversion was not run.}
    \label{tab:fno_training_size}
    {\small\setlength{\tabcolsep}{4pt}%
    \resizebox{\linewidth}{!}{\begin{tabular}{rrcrr@{\hspace{8pt}}rrr@{\hspace{8pt}}rrr}
\toprule
 & & & & & \multicolumn{3}{c}{T2} & \multicolumn{3}{c}{T1} \\
\cmidrule(lr){6-8}\cmidrule(lr){9-11}
Training fields & Width/modes & $\lambda_J$ & $E_F^{\rm val}$ (\%) & $E_J^{\rm val}$ (\%) & $e_m^\star$ (\%) & $k_{\rm or}$ & Final/min & $e_m^\star$ (\%) & $k_{\rm or}$ & Final/min \\
\midrule
1,000 & 8/4 & 10 & 2.646 & 5.450 & 13.91 & 52 & 1.73 & -- & -- & -- \\
10,000 & 8/4 & 10 & 1.357 & 2.493 & 11.79 & 248 & 1.21 & 10.47 & 218 & 1.44 \\
\addlinespace
1,000 & 12/6 & 10 & 1.261 & 1.497 & 11.00 & 465 & 1.38 & -- & -- & -- \\
5,000 & 12/6 & 10 & 1.212 & 1.277 & 10.97 & 371 & 1.33 & -- & -- & -- \\
10,000 & 12/6 & 10 & 0.903 & 0.955 & 10.65 & 551 & 1.23 & 8.85 & 409 & 1.32 \\
\addlinespace
1,000 & 16/8 & 10 & 0.482 & 1.452 & 9.36 & 10000 & 1.00 & -- & -- & -- \\
5,000 & 16/8 & 10 & 0.332 & 0.598 & 9.16 & 10000 & 1.00 & -- & -- & -- \\
10,000 & 16/8 & 10 & 0.264 & 0.387 & 9.60 & 9989 & 1.00 & 7.33 & 9999 & 1.00 \\
\addlinespace
1,000 & 24/10 & 10 & 0.111 & 0.656 & 8.89 & 10000 & 1.00 & 6.62 & 10000 & 1.00 \\
5,000 & 24/10 & 10 & 0.104 & 0.312 & 8.78 & 10000 & 1.00 & 6.54 & 10000 & 1.00 \\
10,000 & 24/10 & 10 & 0.115 & 0.229 & 8.88 & 10000 & 1.00 & 6.56 & 10000 & 1.00 \\
\bottomrule
\end{tabular}
}}
\end{table}

The two smallest models semiconverge even with the Mat\'ern geometry. At
width and modes $8/4$, with $E_F=1.36\%$, the error reaches its minimum of
$11.79\%$ at iteration $248$ and ends at $14.29\%$. At $12/6$ the minimum is
reached at iteration $551$ and the final error is $23\%$ larger. Without
preconditioning, the same models reach their minima after fewer than $120$
iterations and end at $1.8$ to $2.1$ times the minimum. The behaviour
disappears as the surrogate becomes more accurate. At $16/8$ only a slight
deterioration remains without preconditioning, and from $24/10$ onwards, with
$E_F\leq0.12\%$, every minimum lies at the end of the budget.
Reducing the training data has the same effect
(Table~\ref{tab:fno_training_size}). At $8/4$, training on $1000$ instead of
$10\,000$ fields doubles $E_F$ to $2.65\%$, moves the minimum from iteration
$248$ to $52$ and raises the final-to-minimum ratio from $1.21$ to $1.73$. At
$12/6$ the ratio rises from $1.23$ to $1.38$. For the larger models the
iteration remains stable even with $1000$ fields, and at $24/10$ the error is
insensitive to the training size. The deterioration is not specific to the
structured truth, since the small models trained on $10\,000$ fields also
semiconverge on T1, with ratios of $1.44$ and $1.32$.
Semiconvergence without noise is therefore produced by surrogate error alone,
as the level $\Deff$ predicts. The preconditioner delays and reduces it but
does not remove it, and it lowers the final error of the small models but not
their best error. Smaller models and less data degrade the forward and derivative
accuracy together, so the study cannot tell whether $\epsop$ or $\varepsilon_{D,C}$
drives the semiconvergence; each model is also trained only once.

\subsection{Preconditioners satisfying the regularity condition}
\label{sec:reaction-low-smoothness}

The baseline geometry lies outside Theorem~\ref{thm:convergence},
because neither truth satisfies $m^\dagger-m_0\in\M_C$. For T1 the condition
holds when $\nu_C<1/2$, since $\M_C=H^{\nu_C+1}(\Omega)$ and a Mat\'ern
sample with smoothness $3/2$ lies in $H^s(\Omega)$ only for $s<3/2$. T2
remains outside $\M_C$ for every $\nu_C$, since a coefficient with jumps
across curves does not belong to $H^s(\Omega)$ for $s\geq1/2$. We repeat the
noiseless experiments of Section~\ref{sec:reaction-noiseless} with
$\nu_C\in\{0.25,\,0.4\}$ and $\ell_C\in\{0.1,\,0.5,\,1.0\}$ in the kernel \eqref{eq:matern-kernel}. The trained networks, data,
initial field, step rule \eqref{eq:reaction-adaptive-step} and
$10\,000$-iteration budget are unchanged.

\paragraph{Reconstruction accuracy.}
Table~\ref{tab:preconditioner_low_smoothness} compares the smallest
coefficient error reached with $\ell_C=0.1$ with the baseline and with the
identity. With $\ell_C=0.1$, both values of $\nu_C$ improve on the baseline
for every method. For $\nu_C=0.25$, the exact-FEM error decreases from
$4.12\%$ to $3.34\%$ on T1 and from $6.94\%$ to $6.01\%$ on T2, and the
surrogate errors decrease by between $1.5\%$ and $10.5\%$ of their baseline values. Longer length scales, $\ell_C=0.5$ and $1$, generally increase the error
(results not shown). The ordering of
the surrogates on T1 is unchanged, and on T2 only FNO and GenFNO exchange
places. The stability observed in Section~\ref{sec:reaction-geometry} is also
retained. In every one of the $60$ runs the final error is within $5\%$ of the
minimum, whereas without preconditioning GenFNO and GenDeepONet end at up to
three times their minima. The conclusions drawn from the baseline geometry
therefore also hold in geometries covered by the analysis.

\begin{table}[tbp]
    \centering
    \caption{Best relative coefficient error $e_m(k_{\mathrm{or}})$, in percent,
    for noiseless inversion with the baseline Mat\'ern geometry
    ($\nu_C=3/2$), the two geometries with $\nu_C<1/2$ that satisfy the
    regularity condition for T1, and the identity. All Mat\'ern geometries
    use $\ell_C=0.1$. The exact-FEM identity run was performed for T2 only.}
    \label{tab:preconditioner_low_smoothness}
    {\small\begin{tabular}{llrrrr}
\toprule
 & & \multicolumn{3}{c}{Mat\'ern, $\ell_C=0.1$} & \\
\cmidrule(lr){3-5}
Truth & Method & $\nu_C=3/2$ & $\nu_C=0.25$ & $\nu_C=0.4$ & Identity \\
\midrule
T1 & Exact FEM & 4.12 & 3.34 & 2.89 & -- \\
T1 & FNO & 6.64 & 6.05 & 6.16 & 5.32 \\
T1 & DeepONet & 7.65 & 7.50 & 7.51 & 8.77 \\
T1 & GenFNO & 5.83 & 5.74 & 5.70 & 7.93 \\
T1 & GenDeepONet & 7.11 & 6.87 & 6.92 & 8.22 \\
\midrule
T2 & Exact FEM & 6.94 & 6.01 & 5.87 & 1.62 \\
T2 & FNO & 8.86 & 7.93 & 8.16 & 5.63 \\
T2 & DeepONet & 9.70 & 9.39 & 9.44 & 9.57 \\
T2 & GenFNO & 8.57 & 8.33 & 8.34 & 8.80 \\
T2 & GenDeepONet & 9.41 & 9.00 & 9.05 & 9.62 \\
\bottomrule
\end{tabular}
}
\end{table}

\paragraph{Checking the hypotheses along the iteration.}
Since the truth is known, we can evaluate along the trajectories the two
quantities on which the descent argument of Lemma~\ref{lem:descent} rests.
With $e_k=m_k-m_h^\dagger$, the first is the cone quantity of the exact
model, which the tangential cone condition requires to stay below $1/2$,
\[
    \eta_k=
    \frac{\|\F(m_k)-\F(m_h^\dagger)-D\F(m_k)e_k\|_{W_y}}
         {\|\F(m_k)-\F(m_h^\dagger)\|_{W_y}}.
\]
The second measures how well the linearisation of the model
$\F_\star\in\{\F,\F_\theta\}$ that drives the iteration predicts its own
residual $R_k=\F_\star(m_k)-y_h$ in the direction of the error. Together with
the relative step size, it is
\[
    \zeta_k=
    \frac{\|D\F_\star(m_k)e_k-R_k\|_{W_y}}{\|R_k\|_{W_y}},
    \qquad
    s_k=\frac{\alpha_k\,q_k^\top d_k}{\|R_k\|_{W_y}^2}.
\]
Repeating the expansion in the proof of Lemma~\ref{lem:descent} in $\M_C$,
where $\|d_k\|_{\M_C}^2=q_k^\top d_k$, gives
\[
    \|e_{k+1}\|_{\M_C}^2-\|e_k\|_{\M_C}^2
    \leq
    -\alpha_k\|R_k\|_{W_y}^2\,(1-2\zeta_k)
    \qquad\text{whenever } s_k\leq1,
\]
provided neither the $\ell^\infty$ cap nor the projection is active. In
words, a step with $\zeta_k<1/2$ and $s_k\leq1$ reduces the error in
$\M_C$, so $\zeta_k<1/2$ certifies descent at that step. By
Lemma~\ref{lem:linearization}, $\zeta_k$ is at most $\eta$ plus terms
proportional to $1/\|R_k\|$, so the certificate is expected to fail once the
residual approaches the level set by discretisation mismatch and surrogate
error. Both quantities require the truth, so they test the theory and do
not give a stopping rule.

\begin{table}[tbp]
    \centering
    \caption{Diagnostics along the noiseless trajectories with $\nu_C=0.25$
    and $\ell_C=0.1$, evaluated every ten iterations. Here $\zeta_{k_{\mathrm{or}}}$
    is the value of $\zeta_k$ at the oracle iterate and $k_\zeta$ the first
    evaluated iterate with $\zeta_k\geq1/2$ (a dash if there is none). The
    column $e_m(k_\zeta)/e_m(k_{\mathrm{or}})$ is the error at $k_\zeta$ divided
    by the oracle error. The $\M_C$ column states whether $\|e_k\|_{\M_C}$ is
    non-increasing up to the oracle iterate.}
    \label{tab:low_smoothness_diagnostics}
    {\small\begin{tabular}{llrrrrcr}
\toprule
Truth & Method & $\zeta_{k_{\mathrm{or}}}$ & $k_\zeta$ & $k_{\mathrm{or}}$ & $e_m(k_\zeta)/e_m(k_{\mathrm{or}})$ & $\M_C$ monotone & $\max_k s_k$ \\
\midrule
T1 & Exact FEM & 0.68 & 7,330 & 10,000 & 1.04 & yes & 0.38 \\
T1 & FNO & 0.80 & 1,770 & 10,000 & 1.15 & yes & 0.23 \\
T1 & DeepONet & 1.07 & 340 & 10,000 & 1.22 & yes & 0.14 \\
T1 & GenFNO & 1.02 & 80 & 4,826 & 1.76 & no & 0.21 \\
T1 & GenDeepONet & 1.54 & 170 & 7,668 & 1.28 & yes & 0.23 \\
\midrule
T2 & Exact FEM & 0.30 & -- & 10,000 & -- & yes & 0.19 \\
T2 & FNO & 0.85 & 2,260 & 10,000 & 1.13 & yes & 0.18 \\
T2 & DeepONet & 1.28 & 150 & 10,000 & 1.21 & yes & 0.18 \\
T2 & GenFNO & 1.20 & 90 & 3,416 & 1.46 & no & 0.16 \\
T2 & GenDeepONet & 3.11 & 90 & 3,538 & 1.35 & yes & 0.17 \\
\bottomrule
\end{tabular}
}
\end{table}

Table~\ref{tab:low_smoothness_diagnostics} gives the diagnostics for
$\nu_C=0.25$ and $\ell_C=0.1$; the other five geometries behave in the same
way. Three observations follow. First, the hypotheses on the exact model and
the step hold: $\eta_k\leq0.171$ on T1 and $\eta_k\leq0.351$ on T2 in all
$60$ runs, the largest values occurring at the common initial field, and
$s_k\leq0.41$ throughout. Second, the certificate is lost much earlier with
a surrogate. For the exact model $\zeta_k$ reaches $1/2$ late or not at all,
whereas for the surrogates it does so in $43$ of the $48$ runs, after fewer
than $700$ iterations for the generalised architectures and latest for FNO.

Third, the certificate is conservative. The coefficient error continues to decrease
well beyond the first crossing, and in the table stopping there would
increase the error by $13\%$--$22\%$ for FNO and DeepONet and by $28\%$--$76\%$
for the generalised architectures. The error measured in $\M_C$, which
Lemma~\ref{lem:descent} predicts to decrease, is non-increasing
up to the oracle iterate in $54$ of the $60$ runs, including all exact-FEM
runs, even on T2, which lies outside the theory. The six exceptions are the
GenFNO runs with $\ell_C=0.1$, shown in the table, and two GenDeepONet runs
on T2 with larger $\ell_C$, the surrogates whose certificate fails earliest.
The $\ell^\infty$ cap is active in at most $30$ early iterations of each run,
and the projection onto $[0,200]$ is active in up to about a third of the
iterations in some T2 runs with DeepONet and GenDeepONet. Neither is covered
by the analysis.

\subsection{Computational cost}
\label{sec:reaction-computational-cost}

To keep the comparison fair, the exact-FEM baseline assembles the
finite-element operator and its algebraic multigrid preconditioner once per
iterate and reuses them for the forward, tangent and adjoint solves.
A $10\,000$-iteration exact inversion took $1.7$--$2.2$ hours on eight CPU
cores, and a surrogate inversion $1.8$--$4.0$ minutes on one NVIDIA L40 GPU,
a reduction by factors of $26$--$34$ for DeepONet, $34$--$44$ for FNO,
$49$--$63$ for GenFNO and $56$--$72$ for GenDeepONet. These are loop times
at equal iteration counts, not at equal accuracy, comparing a CPU and a GPU
implementation without hardware normalisation. Offline, the $15\,000$
training, validation and test fields with their derivative labels required
$15\,000$ forward and $240\,000$ tangent solves, about $1.06$ hours of
serial-equivalent computation, and training the selected models took between
$2.5$ and $27$ hours on one A100 GPU, excluding hyperparameter search
(for FNO and GenDeepONet the time of the run with $\lambda_J=100$ is used,
since the timing records of the selected runs were not retained).

\subsection{Observational noise and empirical stopping}
\label{sec:reaction-noise}

\begin{figure}[!t]
    \centering
        \includegraphics{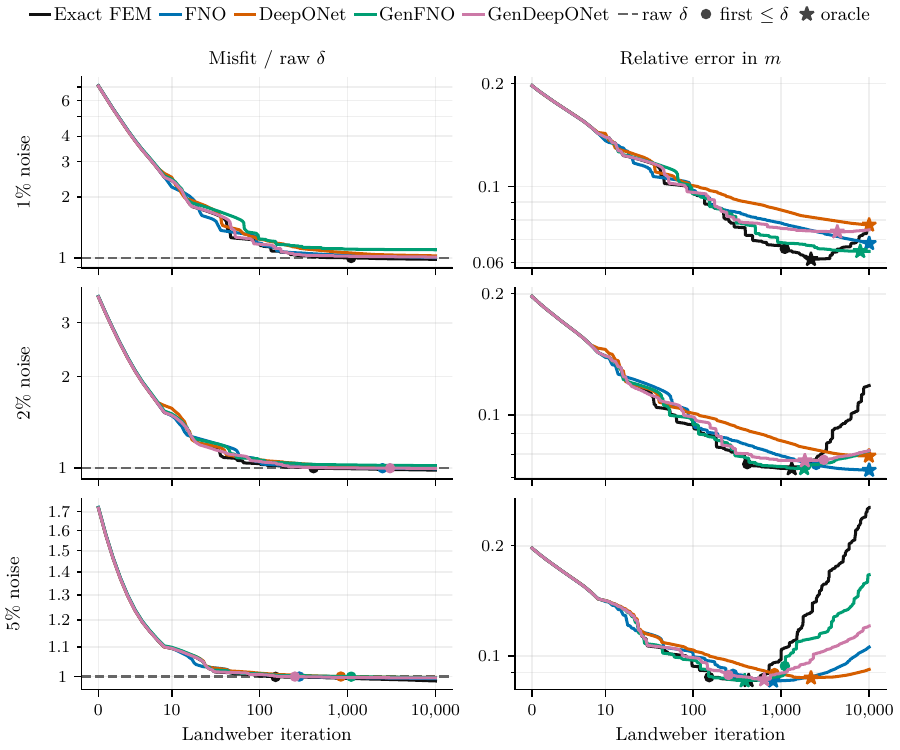}
    
    \caption{Noisy T1 inversions at $1\%$, $2\%$ and $5\%$ observational
    noise. Left: residual divided by the raw noise norm $\delta$; circles mark
    the first raw-$\delta$ crossing. Right: relative coefficient error; stars denote oracle minima. Surrogate curves use the surrogate
    residual and should not be interpreted as FEM re-evaluated physical
    residuals.}
    \label{fig:T1_noisy_convergence}
\end{figure}

We finally examine how observational noise, surrogate error and stopping
interact. For each truth we add one fixed noise realisation at relative
levels $1\%$, $2\%$ and $5\%$, the same for the exact and surrogate
inversions, so these are paired case studies rather than averages over the
noise. Figure~\ref{fig:T1_noisy_convergence} shows the trajectories for T1;
those for T2 behave in the same way.

The analysis stops at the first iterate with $\|R_k\|\leq\tau\Deff$, with
$\tau>2(1+\eta)/(1-2\eta)$. Since $\Deff$ cannot be computed, we apply the
rule retrospectively with the raw level,
$k_\tau=\min\{k:\|R_k\|_{W_y}\leq\tau\delta\}$, for $\tau\in\{1,1.5,2\}$
(Table~\ref{tab:noisy_tau_stopping}). With $\tau=1$ the surrogate iterations
often never stop. At low noise their residuals level off above $\delta$,
since a learned model cannot fit the data more closely than its own error
allows, which is what $\Deff$ accounts for. None of the eight surrogate runs
reaches $\delta$ at $1\%$ noise, although the exact runs do, two of eight do
at $2\%$ and all do at $5\%$. Whenever the rule stops, it stops within $10\%$
of the oracle error. With $\tau=1.5$ or $2$ the rule always stops, but within
the first $100$ iterations, where all models still follow the exact
iteration closely, so every model stops at about the same error, $1.5$ to
$2.3$ times the exact-FEM oracle error. This early stopping reflects the
discrepancy rule on this problem rather than surrogate error. The values of
$\tau$ admitted by the analysis are larger still: for the largest $\eta_k$
measured in Section~\ref{sec:reaction-low-smoothness}, $0.17$ on T1 and
$0.35$ on T2, they exceed $3.5$ and $9$.

\begin{table}[!htb]
    \centering
        \caption{Retrospective discrepancy stopping with the raw noise level.
    Entries are relative coefficient errors, at the oracle iterate
    and at the stopping index $k_\tau$. Surrogate rows give the range over the four surrogates. When
    only some of them reach the threshold, the range is over those that do,
    and their number is given in brackets. For example, (2/4) means that two
    of the four surrogates reach it. A
    dash means that the threshold is not reached within $10\,000$ iterations.
    At $5\%$ noise on T1 with $\tau=2$ the threshold is met by the initial
    field.}
    \label{tab:noisy_tau_stopping}
    \begin{tabular}{lllrrrr}
\toprule
 & & & & \multicolumn{3}{c}{$e_m(k_\tau)$ (\%)} \\
\cmidrule(lr){5-7}
Truth & Noise & Model & $e_m(k_{\mathrm{or}})$ (\%) & $\tau=1$ & $\tau=1.5$ & $\tau=2$ \\
\midrule
T1 & 1\% & Exact FEM & 6.1 & 6.6 & 10.8 & 12.7 \\
 & & Surrogates & 6.5--7.7 & -- & 10.8--11.0 & 12.6--12.9 \\
T1 & 2\% & Exact FEM & 7.4 & 7.6 & 14.0 & 16.1 \\
 & & Surrogates & 7.3--7.9 & 7.5--7.7 (2/4) & 13.7--14.0 & 16.1 \\
T1 & 5\% & Exact FEM & 8.6 & 8.7 & 18.3 & 19.7 \\
 & & Surrogates & 8.5--8.7 & 8.9--9.4 & 18.3 & 19.7 \\
\midrule
T2 & 1\% & Exact FEM & 8.2 & 8.3 & 12.4 & 13.5 \\
 & & Surrogates & 8.8--9.8 & -- & 12.1--12.5 & 13.3--13.5 \\
T2 & 2\% & Exact FEM & 9.0 & 9.1 & 14.5 & 16.6 \\
 & & Surrogates & 9.1--9.9 & -- & 14.1--14.5 & 16.5--16.6 \\
T2 & 5\% & Exact FEM & 10.2 & 10.4 & 18.6 & 20.6 \\
 & & Surrogates & 10.1--10.5 & 10.3--10.8 & 18.6 & 20.6 \\
\bottomrule
\end{tabular}

\end{table}

Beyond the oracle iterate the exact-FEM error grows fastest, more than
doubling by the end of the run at $5\%$ noise on T1, possibly because a
surrogate that resolves fine scales less accurately also fits less of the
noise; with a single noise realisation we note this only as an observation. At $5\%$ noise the oracle errors of the exact and surrogate
iterations differ by less than $0.4$ percentage points on both truths, so
surrogate error need not dominate once the noise is appreciable.

\paragraph{Can held-out surrogate errors be used for stopping?}
The level $\Deff$ involves uniform error bounds, so we tested whether
thresholds built from held-out test errors could replace it on the $18$
noiseless T2 surrogate trajectories that semiconverge. Forward-error
thresholds were missed by some trajectories, even at the maximum test error,
and thresholds combining forward and derivative errors as in the theory were
always crossed, but about $3.3\times10^3$ iterations too early, with errors
$1.4$--$1.8$ times the best attainable; finite-test errors are therefore not
informative about when to stop, which remains an open question in the
noiseless, surrogate-limited regime.

\subsection{Transient Darcy flow}
\label{sec:darcy}

The reaction benchmark compares four architectures on a problem where exact
derivatives are cheap. This example shows that the approach carries over to
a time-dependent problem with sparse space--time data, using a single
standard architecture: a DeepONet, whose trunk can be evaluated at any
space--time point, so that only the observed values need to be computed
during the inversion. The example lies outside the assumptions of
Section~\ref{sec:landweber}, since the surrogate is only piecewise smooth and
the truth is rougher than $\M_C$.

The unknown is a log-permeability $m$ on $\Omega=(0,1)^2$. For
$0<t\le T=0.2$, the pressure $p$ satisfies
\begin{equation}
 \partial_t p-\nabla\cdot\bigl(e^{m}\nabla p\bigr)=f
 \ \text{ in }\Omega\times(0,T],\qquad
 e^{m}\nabla p\cdot n=0\ \text{ on }\partial\Omega\times(0,T],\qquad
 p(\cdot,0)=1 .
 \label{eq:darcy-forward}
\end{equation}
The source $f$ consists of six injectors and six producers of unit rate,
modelled by normalised Gaussians of width $0.04$, so that the net source
vanishes. The pressure is observed at $24$ sensors (the well centres and
twelve further observation wells) at the $51$ times $t_\ell=0.004\,\ell$,
giving $1224$ data, and $\F(m)$ denotes this vector of observations. The
tangent and adjoint equations are given in Appendix~\ref{app:darcy}.

\paragraph{Discretisation and training.}
The parameter lives on a $120\times120$ nodal grid, and the pressure uses
continuous piecewise-linear elements on a $256\times256$-cell mesh with
$400$ backward-Euler steps. Training coefficients are zero-mean Gaussian
Mat\'ern fields with variance $0.25$, $\nu\in\{1,1.5,2\}$ and
$\ell\in\{0.08,\ldots,0.30\}$. We use $20\,000$ training, $2\,000$ validation
and $2\,000$ test fields, each with pressures at the sensors and on a
$65\times65$ grid, and directional derivatives along $8$ of $32$ empirical
KL directions. The surrogate is a standard DeepONet of about
$1.55\times10^{7}$ parameters, with a convolutional branch acting on $m$ and
a trunk acting on space--time coordinates. The branch uses ReLU activations
and max pooling, so the surrogate is only piecewise smooth in $m$. It is
trained with the pressure loss and a directional-derivative loss of weight
$\lambda_J=1$ (not tuned); for comparison, the same architecture is trained
with the pressure loss alone. On the test fields, the median relative forward
error over all sensors and times is $0.61\%$ with derivative training and
$0.66\%$ without it, whereas the median relative error of the directional
derivatives falls from $8.8\%$ to $4.0\%$, and derivative training reduces
it on all but $3$ of the $16\,000$ field--direction pairs.

\paragraph{Inversion setup.}
The truth is a validation field, not used to fit the network weights. Data
are observed on $0\le t\le0.2$ and generated on a refined mesh with
$512\times512$ cells and $800$ steps, so that no method inverts its own
discretisation. Noisy data add a single fixed noise realisation, scaled so
that $\delta=\varrho\,\|y^0-\mathbf 1\|$ for $\varrho\in\{1,2,5\}\%$, where
$y^0$ are the clean observations. Both iterations use the nodal Mat\'ern
preconditioner with $\nu_C=0.5$ and $\ell_C=0.1$, so that
$\M_C=H^{3/2}$, whereas the truth was drawn with $\nu=1$ and $\ell=0.135$
and lies outside it, the step rule \eqref{eq:reaction-adaptive-step} with
each nodal update capped at $0.75$, projection onto $[-2.8,2.8]$ and
$m_0\equiv0$. The exact iteration was run for $1000$ steps and the surrogate
iteration for $20\,000$; with noise we report the first iterate whose own
residual falls below $\delta$ ($\tau=1$), and without noise the last iterate.

The aim is prediction as well as recovery. Each reconstruction is simulated
with the finite-element model up to $t=0.4$, with unchanged rates (S1) and
with a changed schedule (S2) in which, from $t=0.2$, three injectors and
three producers double their rates and the others are shut in. We report the
pooled relative error $\|p_{\mathrm{rec}}-p^\dagger\|/\|p^\dagger-1\|$ of the
pressures at the twelve injection and production wells. The forecasts are
finite-element simulations of the recovered field, since the surrogate was
trained only up to $t=0.2$.

\begin{figure}[tbp]
    \centering
    \includegraphics[width=\textwidth]{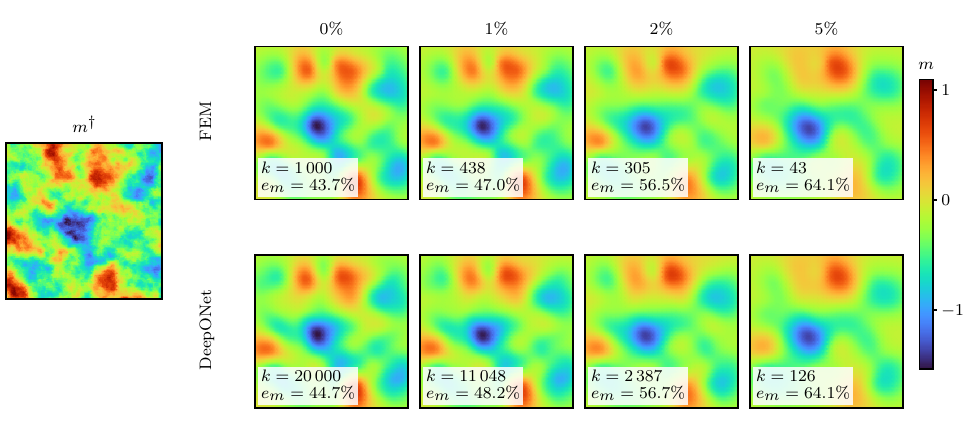}
    \caption{Transient Darcy flow: truth (left) and reconstructions by the
    exact (top) and DeepONet (bottom) iterations at each noise level, on a
    common colour scale, at the iterates of Table~\ref{tab:darcy}.}
    \label{fig:darcy_reconstructions}
\end{figure}

\begin{figure}[tbp]
    \centering
    \includegraphics[width=\textwidth]{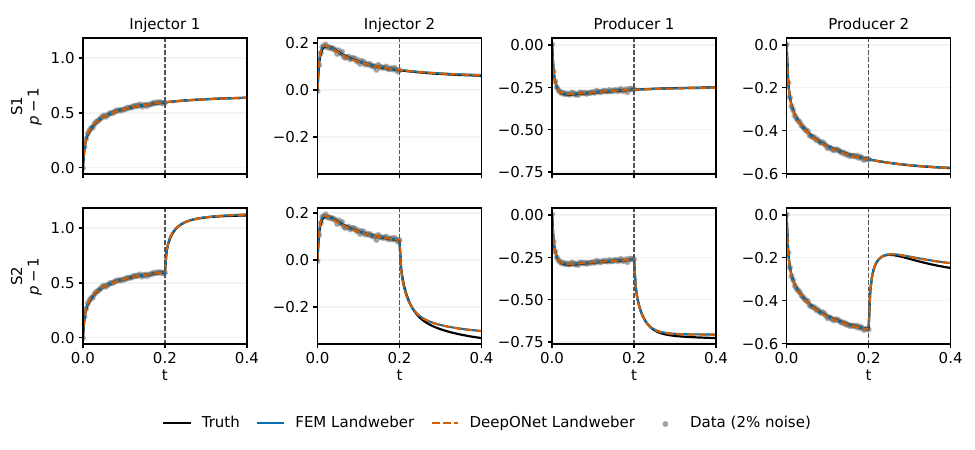}
    \caption{Transient Darcy flow: pressures at two injectors and two
    producers, chosen before the forecasts were computed, from FEM
    simulation of the fields recovered from data with $2\%$ noise. Top: S1,
    unchanged rates. Bottom: S2, where from $t=0.2$ the first injector and
    producer double their rates and the second are shut in. The dashed line
    marks the end of the data window.}
    \label{fig:darcy_wells}
\end{figure}

\paragraph{Results.}
With noisy data the surrogate iteration stops by the discrepancy principle
at every noise level, and its reconstruction error is within $1.3$
percentage points of the exact iteration
(Table~\ref{tab:darcy}); the reconstructions are visually almost
indistinguishable (Figure~\ref{fig:darcy_reconstructions}). The well
pressures agree as closely: both methods reproduce the history to about
$1\%$, forecast S1 to below $1\%$ and the changed schedule S2 to
$1.3$--$3.2\%$, with no systematic advantage for either
(Figure~\ref{fig:darcy_wells}). The error in the recovered field is large,
between $44\%$ and $64\%$, because only $24$ sensors inform a rough field,
but the predictions of interest are accurate, including under a schedule
that the data do not contain.

\begin{table}[tbp]
    \centering
    \caption{Transient Darcy flow: exact FEM and DeepONet-surrogate Landweber
    iteration. With noise, the first iterate whose own residual is below
    $\delta$; without noise, the last iterate. Time is the measured
    loop time (exact on CPU, surrogate on one L40 GPU). $r_k/\delta$ is the
    FEM residual of the stopped iterate relative to the noise level.
    Well-pressure errors are relative errors at the twelve injection and
    production wells, from FEM simulation of the recovered field, over the
    history window and for the forecasts S1 (unchanged rates) and S2
    (changed rates) on $0.2<t\le0.4$.}
    \label{tab:darcy}
    \begin{tabular}{llrrrrrrr}
\toprule
 & & & & & & \multicolumn{3}{c}{Well-pressure error (\%)} \\
\cmidrule(lr){7-9}
Noise & Method & $k$ & Time (s) & $e_m$ (\%) & $r_k/\delta$ & History & S1 & S2 \\
\midrule
0\% & FEM  & 1000  & 22997 & 43.75 & --    & 0.23 & 0.15 & 1.43 \\
    & DeepONet & 20000 & 293.5 & 44.74 & --    & 0.43 & 0.39 & 1.29 \\
\addlinespace
1\% & FEM  & 438   & 9895  & 46.96 & 0.996 & 0.26 & 0.24 & 2.17 \\
    & DeepONet & 11048 & 160.9 & 48.22 & 1.075 & 0.45 & 0.35 & 1.97 \\
\addlinespace
2\% & FEM  & 305   & 6845  & 56.54 & 1.000 & 0.48 & 0.58 & 2.86 \\
    & DeepONet & 2387  & 34.9  & 56.68 & 1.020 & 0.59 & 0.50 & 2.91 \\
\addlinespace
5\% & FEM  & 43    & 982.9 & 64.13 & 1.000 & 1.05 & 0.97 & 3.10 \\
    & DeepONet & 126   & 2.66  & 64.08 & 1.004 & 1.08 & 0.97 & 3.20 \\
\bottomrule
\end{tabular}

\end{table}

It is worth noting three further points. First, the surrogate stops on its
own residual, and
the FEM residual of the stopped iterate exceeds $\delta$ by $7\%$, $2\%$ and
$0.4\%$ at the three noise levels, so the discrepancy principle holds only
approximately in the physical model. Second, the surrogate iteration needs
$25$, $8$ and $3$ times as many steps as the exact one, consistent with the
level $\Deff$, which makes the iteration slow down as $\delta$ approaches the
surrogate error. Third, without noise the surrogates trained with and
without the derivative loss reach $44.7\%$ and $44.9\%$, against $43.7\%$
for the exact iteration, so derivative training, which halves the derivative
error, does not change the reconstruction, as in the reaction benchmark.

\paragraph{Cost.}
A surrogate step takes about $15$\,ms on the GPU, against about $23$\,s for
an exact step (two forward, one tangent and one adjoint solve over $400$
time steps). Despite the larger number of steps, the time to stopping falls
from $2.7$\,h, $1.9$\,h and $16$\,min to $2.7$\,min, $35$\,s and $3$\,s at
$1\%$, $2\%$ and $5\%$ noise; the exact solver ran serially, so these ratios
compare two implementations, not two optimised codes. Offline, generating the
$24\,000$ fields with derivative labels took $1453$ allocated CPU-core hours,
as much as about $66$, $96$ and $660$ exact inversions, and training took
$19$ GPU hours. The surrogate therefore pays off when the inverse problem is
solved repeatedly with different data from the same wells, and, since this
cost is paid before the data arrive, when the field and its predictions must
be updated as new measurements come in. The surrogate can also be queried
away from the sensors without retraining: on $200$ validation fields its
median pressure error on the $65\times65$ training grid is $0.91\%$, against
$0.62\%$ at the sensors, and a full-grid evaluation takes $0.16$\,s against
$3.4$\,s for a finite-element solve, so pressures at, for instance, planned
observation wells are available at a small fraction of the cost (we have
checked the accuracy only on the training grid).

\FloatBarrier
\section{Conclusions}
\label{sec:conclusions}

We have studied nonlinear Landweber iteration in which the forward map and
its derivative are replaced by a differentiable neural-operator surrogate.
Learning the inverse map directly is not a viable alternative: near a
discontinuity, no continuous network approximates it uniformly
(Proposition~\ref{prop:inverse-blowup}). The forward map, by contrast, can be
learned together with its derivative by encoder--network--decoder operators,
a class that generalised DeepONets and fully coupled Fourier operators
realise exactly. For such a surrogate, a discrepancy principle with the level
$\Deff$, which combines the noise level with the forward and derivative
errors, gives descent before stopping and finite termination. With a compact
preconditioner, the stopped reconstructions converge as the data and
surrogate errors vanish, provided $m^\dagger-m_0\in\M_C$
(Theorem~\ref{thm:convergence}).

The experiments support the central mechanism of the analysis: surrogate
error acts like data noise. Without preconditioning, three of the four
surrogates in the reaction benchmark semiconverge early even with exact data,
and develop spurious fine-scale structure, whereas the exact iteration does
not. The compact preconditioner removes this instability by damping the
poorly learned components of the adjoint direction, but it also slows the
recovery of sharp features. Derivative-informed training along a fixed set
of directions reduced the derivative test error by $45$--$94\%$ but hardly
changed the best reconstructions, in both the reaction and the Darcy
examples; its clearer effect was to delay semiconvergence. What matters for
inversion is therefore not the derivative accuracy along the training
directions but the accuracy of the learned adjoint direction in the geometry
of the iteration. Once the offline cost is paid, the surrogate iterations
were between about $25$ and $300$ times faster than the finite-element ones,
which pays off when the inverse problem is solved repeatedly or must be
updated as new data arrive.

Several questions remain open. The level $\Deff$ involves uniform error
bounds that are not computable in practice, and thresholds built from
held-out test errors did not identify when to stop in the noiseless,
surrogate-limited regime. A posteriori estimates of the surrogate error along
the inverse trajectory, possibly combined with occasional finite-element
evaluations, are a natural next step. A second direction is training that
targets the weighted derivative error $(D\F_\theta-D\F)C^{1/2}$ directly, so
that the training and inversion geometries are aligned. The analysis is
carried out in Hilbert spaces; its extension to Banach spaces along the lines
of \citet{SKHK2012} would allow penalties better suited to sharp interfaces.
Finally, although a discontinuous inverse map cannot be learned directly, a
regularised neural operator that balances the approximation and
noise-propagation errors may be feasible.

\section{Declarations}

Throughout this work, generative AI was used to polish the wording and
correct grammar, to proofread parts of the manuscript and suggest
clarifications and modifications, to draft figure captions and parts of the
descriptions of the numerical experiments, and to write all software used in
this paper, following our specifications, including the finite-element
solvers, the neural architectures, the training code and the Landweber
implementations.
The tools were ChatGPT (OpenAI; Astra) and OpenAI Codex (based on GPT-6; desktop application version 26.1016.62342) for
the software, and Claude (Anthropic; Opus~5.5) for the text. No AI system ran
jobs on the computing cluster, and every reported number was computed by
archived scripts from saved results. The solvers were checked by mesh
refinement, finite differences and the adjoint identity.
The research questions, the theoretical framework and results, the choice of
problems and architectures, and all decisions on what to report are our own.
We carefully checked all AI-generated content and suggestions and take full
responsibility for all mathematical content, software, numerical results,
interpretations, and conclusions.

\paragraph{Code and data availability.}
The code for the reaction benchmark, including scripts that regenerate all
its tables and figures, is at
\url{https://github.com/Marco-Iglesias-Nottingham/neural-operator-landweber-reaction},
with data, trained models and saved results at
\url{https://doi.org/10.5281/zenodo.23767592}. The code and data for the
Darcy example are at
\url{https://github.com/Marco-Iglesias-Nottingham/neural-operator-landweber-darcy}
and \url{https://doi.org/10.5281/zenodo.23879298}.

\acks{K. Gai and S. Lu were supported by National Key Research and Development Programs of China (No. 2023YFA1009103), NSFC (No. 12631015 and No. 92570106).
Simulations were performed using the Sulis Tier~2 HPC platform hosted
by the Scientific Computing Research Technology Platform at the University of
Warwick. Sulis is funded by EPSRC, United Kingdom, Grant EP/T022108/1 and the
HPC Midlands+ consortium.
}

\appendix
\section{Approximation proofs}
\label{app:approx-proofs}

This appendix gives full proofs of Theorems~\ref{thm:compact-directions}
and~\ref{thm:operator-approx}. Both proofs follow the same four steps:
finite-rank truncation of input and output, control of the derivative of the
truncated operator, reduction to a map between Euclidean spaces, and
approximation of that map by a network.

\subsection{Background results}

\begin{definition}[Approximation property]
\label{def:ap}
A Banach space $X$ has the approximation property if, for every compact set
$K\subset X$ and every $\epsilon>0$, there is a bounded finite-rank operator
$T:X\to X$ with $\sup_{x\in K}\|x-Tx\|_X\le\epsilon$.
\end{definition}

\begin{lemma}[{\citealp[Lemma~18]{Kovachki2023}}]
\label{lem:schauder-ap}
Every Banach space with a Schauder basis has the approximation property.
More precisely, if $P_n$ are the partial-sum projections of the basis, then
$L:=\sup_n\|P_n\|<\infty$ and $P_n\to I$ uniformly on every compact set.
\end{lemma}

\begin{lemma}[{\citealp[Lemma~21]{Kovachki2023}}]
\label{lem:compact-union}
Let $X,Y$ be Banach spaces, $K\subset X$ compact, and $T,T_n:X\to Y$
continuous with $\sup_{x\in K}\|T(x)-T_n(x)\|_Y\to0$. Then
$T(K)\cup\bigcup_{n\ge1}T_n(K)$ is compact in $Y$.
\end{lemma}

\begin{lemma}[Finite-dimensional $C^1$ approximation]
\label{lem:finite-c1}
Let $\sigma\in C^1(\mathbb R)$ be nonpolynomial, let $V\subset\mathbb R^r$ be
open, let $g\in C^1(V,\mathbb R^s)$, and let $B\subset V$ be compact. For every
$\epsilon>0$ there is a one-hidden-layer network
$f_\theta(z)=W_2\,\sigma(W_1z+b_1)+b_2$, with $\sigma$ applied componentwise,
such that $\sup_{z\in B}\|g(z)-f_\theta(z)\|<\epsilon$ and
$\sup_{z\in B}\|Dg(z)-Df_\theta(z)\|<\epsilon$.
\end{lemma}
\begin{proof}
Let $\chi\in C_c^\infty(V)$ equal one on a neighbourhood of $B$. Extended by
zero, $\chi g$ belongs to $C^1(\mathbb R^r,\mathbb R^s)$ and coincides with $g$
near $B$. By \citet[Theorem~4.1]{Pinkus1999}, for $\sigma\in C^1$ nonpolynomial,
linear combinations of ridge functions $\sigma(w\cdot z+b)$ are dense in
$C^1(\mathbb R^r)$ uniformly on compact sets. Apply this to each component of
$\chi g$ with a tolerance small enough that the Euclidean and Frobenius norms
of the value and Jacobian errors are below $\epsilon$ on $B$, and stack the $s$
hidden layers into one.
\end{proof}

\subsection{Proof of Theorem~\ref{thm:compact-directions}}
\phantomsection
\label{proof:compact-directions}

Throughout, we let $\epsilon>0$ be arbitrary and write $D(m):=D\F(m)$.

\paragraph{Step 1: finite-rank truncations and compact sets.}
By the approximation property of $\M$, applied with tolerances $1/n$ on the
compact set $\Adm\cup H$, there are bounded finite-rank operators
$P_n:\M\to\M$ with $\rho_n:=\sup_{m\in\Adm}\|m-P_nm\|_\M\to0$ and
$\rho_n^H:=\sup_{h\in H}\|h-P_nh\|_\M\to0$.
Since $\Adm$ is compact and $U$ open, $\dist(\Adm,\M\setminus U)>0$; discarding
finitely many $n$, we may assume $P_n\Adm\subset U$ for all $n$. By
Lemma~\ref{lem:compact-union}, the sets $Z:=\Adm\cup\bigcup_{n}P_n(\Adm)\subset U$ and
$Z_H:=H\cup\bigcup_nP_n(H)$ are compact. Hence $\F$ and $D$ are uniformly continuous on $Z$, with moduli
of continuity $\omega_\F$ and $\omega_D$, and $L_D:=\sup_{a\in Z}\|D(a)\|<\infty$,
$R_H:=\max\{1,\sup_{h\in Z_H}\|h\|_\M\}<\infty$.
The map $(a,h)\mapsto D(a)h$ is continuous on $Z\times Z_H$, because
$\|D(a)h-D(a')h'\|\le\|D(a)-D(a')\|\,\|h\|+\|D(a')\|\,\|h-h'\|$. Therefore
$W:=\F(Z)\cup\{D(a)h:\ a\in Z,\ h\in Z_H\}$ is compact in $\Y$. By the approximation property of $\Y$, choose a bounded
finite-rank operator $Q:\Y\to\Y$ with
\begin{equation}
 \sup_{y\in W}\|y-Qy\|_\Y\le\epsilon/30 .
 \label{eq:app-Q-tol}
\end{equation}
Since $Z$ and $Z_H$ do not depend on the index $n$, neither does $W$, and $n$
may be chosen after $Q$. Choose $n$ so large that
\begin{equation}
 \omega_\F(\rho_n)\le\epsilon/6,\qquad
 \omega_D(\rho_n)R_H\le\epsilon/12,\qquad
 L_D\,\rho_n^H\le\epsilon/12,
 \label{eq:app-n-tol}
\end{equation}
write $P=P_n$, and set $\F_{P,Q}:=Q\circ\F\circ P$ on $\{m:Pm\in U\}$.

\paragraph{Step 2: accuracy of the truncated operator.}
For $m\in\Adm$, since $Pm\in Z$,
\[
 \|\F(m)-\F_{P,Q}(m)\|_\Y
 \le\|\F(m)-\F(Pm)\|_\Y+\|(I-Q)\F(Pm)\|_\Y
 \le\omega_\F(\rho_n)+\epsilon/30<\epsilon/3 .
\]
By the chain rule, $D\F_{P,Q}(m)=Q\,D(Pm)\,P$, and
\begin{equation}
 D(m)-Q\,D(Pm)\,P
 =(I-Q)D(m)+Q\bigl[D(m)-D(Pm)\bigr]P+Q\,D(m)(I-P).
 \label{eq:app-three-terms}
\end{equation}
Fix $m\in\Adm$ and $h\in H$. We bound the three terms applied to $h$ without
using any bound on $\|Q\|$: in each case we add and subtract terms so that
$I-Q$ acts only on elements of $W$.

First, $D(m)h\in W$, so $\|(I-Q)D(m)h\|_\Y\le\epsilon/30$.

Second, since $m,Pm\in Z$ and $Ph\in Z_H$, both $D(m)Ph$ and $D(Pm)Ph$ lie in
$W$, and
\begin{align*}
 \bigl\|Q[D(m)-D(Pm)]Ph\bigr\|_\Y
 &\le\|(Q-I)D(m)Ph\|_\Y+\|[D(m)-D(Pm)]Ph\|_\Y\\
 &\qquad+\|(I-Q)D(Pm)Ph\|_\Y \le\tfrac{2\epsilon}{30}+\omega_D(\rho_n)R_H .
\end{align*}

Third, $D(m)h$ and $D(m)Ph$ lie in $W$, and
\begin{align*}
 \|Q\,D(m)(h-Ph)\|_\Y
 &\le\|(Q-I)D(m)h\|_\Y+\|D(m)(h-Ph)\|_\Y+\|(I-Q)D(m)Ph\|_\Y\\
 &\le\tfrac{2\epsilon}{30}+L_D\,\rho_n^H .
\end{align*}
Summing and using \eqref{eq:app-n-tol} gives
$\sup_{m\in\Adm,\,h\in H}\|(D(m)-D\F_{P,Q}(m))h\|_\Y
\le5\epsilon/30+\epsilon/12+\epsilon/12=\epsilon/3$.

\paragraph{Step 3: reduction to Euclidean spaces.}
Choose a basis $e_1,\ldots,e_r$ of $\ran P$ and let $S:\mathbb R^r\to\ran P$,
$Sc=\sum_ic_ie_i$. Then $\Psi:=S^{-1}P\in\Lop(\M,\mathbb R^r)$ and $P=S\Psi$.
Similarly, choose a basis $b_1,\ldots,b_s$ of $\ran Q$, let
$\Phi c=\sum_jc_jb_j$, and let $\Lambda:=\Phi^{-1}Q\in\Lop(\Y,\mathbb R^s)$, so
that $Q=\Phi\Lambda$. The set $V:=\{z\in\mathbb R^r:Sz\in U\}$ is open, and
$g:V\to\mathbb R^s$, $g(z):=\Lambda\,\F(Sz)$, belongs to $C^1$ with
$Dg(z)=\Lambda\,D(Sz)\,S$. By construction,
$\F_{P,Q}=\Phi\circ g\circ\Psi$ with $D\F_{P,Q}(m)=\Phi\,Dg(\Psi m)\,\Psi$,
and $B:=\Psi(\Adm)$ is a compact subset of $V$, since $S\Psi\Adm=P\Adm\subset U$.

\paragraph{Step 4: network approximation.}
Let $M_\Psi:=\max\{1,\sup_{h\in H}\|\Psi h\|\}$ and
$M_\Phi:=\max\{1,\|\Phi\|\}$. Lemma~\ref{lem:finite-c1} gives a network
$f_\theta$ whose value and Jacobian errors on $B$ are below
$\delta':=\epsilon/(3M_\Phi M_\Psi)$.
Set $\F_\theta:=\Phi\circ f_\theta\circ\Psi$, an operator of the form
\eqref{eq:rb}. For $m\in\Adm$ and $h\in H$,
$\|\F_{P,Q}(m)-\F_\theta(m)\|_\Y\le\|\Phi\|\delta'<\epsilon/3$ and
$\|(D\F_{P,Q}(m)-D\F_\theta(m))h\|_\Y\le\|\Phi\|\,\delta'\,\|\Psi h\|<\epsilon/3$.
Combining with Step~2 gives both conclusions of
Theorem~\ref{thm:compact-directions}, with errors below $2\epsilon/3$.
\hfill$\square$

\subsection{Proof of Theorem~\ref{thm:operator-approx}}
\label{proof:operator-approx}

Again let $\epsilon>0$ and write $D(m):=D\F(m)$.

\paragraph{Step 1: basis projections.}
Let $(e_j)$ be a Schauder basis of $\M$ with coordinate functionals
$(e_j^*)\subset\M^*$, and $(u_j)$ a Schauder basis of $\Y$ with coordinate
functionals $(u_j^*)$. Let
$P_nx=\sum_{j=1}^ne_j^*(x)e_j$ and $Q_Ny=\sum_{j=1}^Nu_j^*(y)u_j$ be the partial-sum projections, and let $L\ge1$ bound $\|P_n\|$ and $\|Q_N\|$
for all $n,N$ (Lemma~\ref{lem:schauder-ap}). Then $P_n\to I$ uniformly on
compact subsets of $\M$, and $Q_N\to I$ uniformly on compact subsets of $\Y$.
Moreover $P_n^*f=\sum_{j\le n}f(e_j)e_j^*$ for $f\in\M^*$. Since $\M$ is
reflexive, every Schauder basis of $\M$ is shrinking \citep{AlbiacKalton2006},
so $(e_j^*)$ is a basis of $\M^*$ and $P_n^*\to I$ strongly on $\M^*$, with
$\|P_n^*\|=\|P_n\|\le L$.

As in the previous proof, write $\rho_n=\sup_{m\in\Adm}\|m-P_nm\|$, discard
finitely many $n$ so that $P_n\Adm\subset U$, and let
$Z=\Adm\cup\bigcup_nP_n(\Adm)\subset U$, which is compact by
Lemma~\ref{lem:compact-union}. Let $\omega_\F,\omega_D$ be moduli of
continuity of $\F$ and $D$ on $Z$. For $m\in\Adm$,
\begin{equation}
 \|\F(m)-Q_N\F(P_nm)\|_\Y\le\omega_\F(\rho_n)+\sup_{y\in\F(Z)}\|y-Q_Ny\|_\Y .
 \label{eq:app-op-values}
\end{equation}

\paragraph{Step 2: the derivative in operator norm.}
The chain rule and the decomposition \eqref{eq:app-three-terms}, with
$P=P_n$ and $Q=Q_N$, give
$\|D(m)-Q_ND(P_nm)P_n\|\le T_1+T_2+T_3$, where
$T_1=\|Q_N[D(m)-D(P_nm)]P_n\|$, $T_2=\|(I-Q_N)D(m)\|$ and
$T_3=\|Q_ND(m)(I-P_n)\|$.
\emph{First term.} $T_1\le L^2\,\omega_D(\rho_n)$ for all $m\in\Adm$.

\emph{Second term.} Let $B_\M$ be the closed unit ball of $\M$. The set
$\{D(m):m\in\Adm\}$ is compact in $\Lop(\M,\Y)$; for $\eta>0$ choose a finite
$\eta$-net $D(m_1),\ldots,D(m_t)$ with $m_i\in\Adm$. Each $D(m_i)B_\M$ is
relatively compact because $D(m_i)$ is compact, and every $D(m)x$ with
$m\in\Adm$, $x\in B_\M$ lies within $\eta$ of $\bigcup_iD(m_i)B_\M$. Hence
$S:=\bigcup_{m\in\Adm}D(m)B_\M$ is totally bounded, so $\overline S$ is
compact, and
$\sup_{m\in\Adm}T_2=\sup_{y\in S}\|(I-Q_N)y\|_\Y\to0$ as $N\to\infty$.

\emph{Third term.} Since $\|A\|=\|A^*\|$,
$T_3\le L\,\|D(m)(I-P_n)\|=L\,\|(I-P_n^*)D(m)^*\|$. By Schauder's theorem,
$D(m)^*$ is compact, so $D(m)^*B_{\Y^*}$ is relatively compact, and the
uniformly bounded operators $P_n^*$ converge to $I$ uniformly on it. Hence
$\|(I-P_n^*)D(m)^*\|\to0$ for each $m\in\Adm$. The set
$\{D(m)^*:m\in\Adm\}$ is compact in $\Lop(\Y^*,\M^*)$; if
$D(m_1)^*,\ldots,D(m_t)^*$ is an $\eta$-net, then
$\|(I-P_n^*)D(m)^*\|\le\max_i\|(I-P_n^*)D(m_i)^*\|+(1+L)\eta$, so $\sup_{m\in\Adm}T_3\to0$ as $n\to\infty$.

\emph{Choice of $N$ and $n$.} The sets $\F(Z)$ and $\overline S$ do not depend
on $n$. First choose $N$ such that $\sup_{y\in\F(Z)}\|y-Q_Ny\|\le\epsilon/4$
and $\sup_\Adm T_2\le\epsilon/6$. Then choose $n$ such that
$\omega_\F(\rho_n)\le\epsilon/4$, $L^2\omega_D(\rho_n)\le\epsilon/6$ and
$\sup_\Adm T_3\le\epsilon/6$. With $\F_{n,N}:=Q_N\circ\F\circ P_n$,
\eqref{eq:app-op-values} and the three bounds give
$\sup_{m\in\Adm}\|\F(m)-\F_{n,N}(m)\|_\Y\le\epsilon/2$ and
$\sup_{m\in\Adm}\|D(m)-D\F_{n,N}(m)\|_{\Lop(\M,\Y)}\le\epsilon/2$.

\paragraph{Step 3: reduction to Euclidean spaces.}
Let $\Psi:=(e_1^*,\ldots,e_n^*):\M\to\mathbb R^n$, $Sc:=\sum_{j\le n}c_je_j$,
$\Lambda:=(u_1^*,\ldots,u_N^*):\Y\to\mathbb R^N$ and
$\Phi c:=\sum_{j\le N}c_ju_j$, so that $P_n=S\Psi$ and $Q_N=\Phi\Lambda$. On
the open set $V:=\{z\in\mathbb R^n:Sz\in U\}$, the map $g(z):=\Lambda\F(Sz)$ is
$C^1$ by the chain rule, with
$\F_{n,N}=\Phi\circ g\circ\Psi$ and $D\F_{n,N}(m)=\Phi\,Dg(\Psi m)\,\Psi$.
The set $B:=\Psi(\Adm)\subset V$ is compact.

\paragraph{Step 4: network approximation.}
Lemma~\ref{lem:finite-c1} gives a network $f_\theta$ with value and Jacobian
errors below
$\delta':=\epsilon/(2\max\{1,\|\Phi\|\}\max\{1,\|\Psi\|\})$ on $B$. For
$\F_\theta:=\Phi\circ f_\theta\circ\Psi$ and $m\in\Adm$,
$\|\F_{n,N}(m)-\F_\theta(m)\|_\Y\le\|\Phi\|\delta'<\epsilon/2$ and
$\|D\F_{n,N}(m)-D\F_\theta(m)\|\le\|\Phi\|\,\delta'\,\|\Psi\|<\epsilon/2$.
Together with Step~2, both errors are below $\epsilon$.
\hfill$\square$

\section{Standard FNO calculus and Lipschitz estimates}
\label{app:fno-calculus}

Let $\F_\theta^{\mathrm{FNO}}:m\mapsto u$ denote the full standard FNO
map, with $\theta$ collecting its trainable parameters. Consider bounded
linear lifting and output maps $P,Q$, and layers indexed by
$t=0,\ldots,T-1$: $v_0=Pm+p$, $z_t=(W_t+\mathcal K_t)v_t+b_t$ and
$v_{t+1}=\sigma(z_t)$, with output $u=\F_\theta^{\mathrm{FNO}}(m)=Qv_T(m)+q$.
The derivative $D\F_\theta^{\mathrm{FNO}}(m)$ is taken with respect to the
input $m$, with $\theta$ fixed. On a finite grid with a smooth activation,
the chain rule gives $\delta v_{t+1}=D_t(W_t+\mathcal K_t)\delta v_t$
with $D_t w=\sigma'(z_t)w$, and therefore
\begin{equation}
 D\F_\theta^{\mathrm{FNO}}(m)
 =Q D_{T-1}(W_{T-1}+\mathcal K_{T-1})\cdots
       D_0(W_0+\mathcal K_0)P.
 \label{eq:fno-derivative}
\end{equation}
For nonlinear lifting or decoding, their derivatives replace $P,Q$ at
the corresponding intermediate values. The same formula holds on function
spaces whenever these compositions are Fr\'echet differentiable there,
for example with appropriate smooth composition maps on a Sobolev algebra.
Scalar smoothness and an $L^2$ mapping property alone do not establish
Fr\'echet differentiability of the Nemytskii map $L^2\to L^2$.

On $L^2$ of the torus, or a grid with a unitary discrete Fourier transform,
a truncated multiplier satisfies
$\|\mathcal K_t v\|^2=\sum_{k\in\Lambda_t}\|R_t(k)\widehat v(k)\|_2^2
\le\max_{k\in\Lambda_t}\|R_t(k)\|_2^2\|v\|^2$. Also $\|v\mapsto W_tv\|\le\|W_t\|_2$. If the activation is globally
$L_\sigma$-Lipschitz, repeated composition consequently gives
\begin{equation}
 \Lip(\F_\theta^{\mathrm{FNO}})
 \le L_Q L_P L_\sigma^T
       \prod_{t=0}^{T-1}
       \left(\|W_t\|_2+\max_{k\in\Lambda_t}\|R_t(k)\|_2\right).
 \label{eq:fno-lipschitz}
\end{equation}
Here $L_P,L_Q$ are the lifting and decoding Lipschitz constants; for an
$M$-layer perceptron, $L_P\le L_\sigma^{M-1}\prod_{j=1}^M\|P_j\|_2$,
and similarly for $Q$. Biases cancel in layerwise differences. These
are upper bounds, not estimates of the realised Lipschitz constant.
They apply in the indicated $L^2$ or grid norms; a different norm requires
the corresponding boundedness estimates.

\section{Tangent and adjoint for the Darcy example}
\label{app:darcy}

For a perturbation $h$, the tangent $w=Dp(m)h$ of \eqref{eq:darcy-forward}
satisfies
\[
 \partial_t w-\nabla\cdot\bigl(e^{m}\nabla w\bigr)
 =\nabla\cdot\bigl(e^{m}h\nabla p\bigr),\qquad
 \bigl(e^{m}\nabla w+e^{m}h\nabla p\bigr)\cdot n=0,\qquad w(\cdot,0)=0,
\]
and $(D\F(m)h)_{\ell j}=w(s_j,t_\ell)$. For a data-space vector
$r=(r_{\ell j})$, the adjoint $z$ runs backwards in time,
\[
 -\partial_t z-\nabla\cdot\bigl(e^{m}\nabla z\bigr)
 =\sum_{\ell:\,t_\ell>0}\sum_{j}r_{\ell j}\,\delta_{s_j}\,\delta_{t_\ell},
 \qquad e^{m}\nabla z\cdot n=0,\qquad z(\cdot,T^+)=0,
\]
so that each observation time contributes a jump of $z$. With the
$L^2(\Omega)$ inner product on the parameter,
\[
 D\F(m)^*r=-\int_0^T e^{m}\,\nabla p\cdot\nabla z\,\mathrm dt .
\]
The observation at $t=0$ does not depend on $m$. Point evaluation is not
bounded on $H^1(\Omega)$ in two dimensions, so these formulas are formal. The
computations use the discrete tangent and adjoint of the backward-Euler
finite-element scheme, for which the corresponding identities hold exactly.
The discrete adjoint reuses the forward LU factorisation, and the adjoint
action is accumulated on the parameter grid through the transpose of the
interpolation.

\setlength{\bibsep}{1.5pt plus 0.3ex}
\bibliography{references}
\end{document}